\documentclass[preprint]{elsarticle}

\usepackage{amsmath, amssymb, amsthm}
\usepackage{color}
\usepackage{latexsym}
\usepackage{indentfirst}
\usepackage{graphicx}
\usepackage{placeins}
\usepackage{booktabs}
\usepackage{algorithm}
\usepackage{algorithmic}
\usepackage{multirow}
\usepackage{subfig}

\biboptions{compress}

\theoremstyle{plain}

\theoremstyle{definition}

\theoremstyle{remark}
\newtheorem*{remark}{Remark}

\DeclareMathOperator{\diag}{diag}
\DeclareMathOperator{\spanop}{span}

\newcommand{\ud}{\,\mathrm{d}}

\newcommand{\bd}[1]{\boldsymbol{#1}}
\newcommand{\wt}[1]{\widetilde{#1}}

\newcommand{\mc}[1]{\mathcal{#1}}

\newcommand{\abs}[1]{\left\lvert#1\right\rvert}
\newcommand{\norm}[1]{\left\lVert#1\right\rVert}
\newcommand{\average}[1]{\left\langle#1\right\rangle}
\newcommand{\bra}[1]{\langle#1\rvert}
\newcommand{\ket}[1]{\lvert#1\rangle}

\newcommand{\jump}[1]{\big[\hspace{-0.7mm} \big[ #1 \big]
  \hspace{-0.7mm} \big]} 
\newcommand{\mean}[1] {\big\{ \hspace{-0.7mm} \big\{ #1 \big\}
  \hspace{-0.7mm} \big\}}

\newcommand{\KS}{\mathrm{KS}}
\newcommand{\exc}{\epsilon_{\mathrm{xc}}}
\newcommand{\ext}{\mathrm{ext}}
\newcommand{\eff}{\mathrm{eff}}
\newcommand{\DG}{\mathrm{DG}}

\newcommand{\barint}{\kern4pt \raise3.4pt\hbox{\vrule height.6pt
    width7pt} \kern-11pt \int}

\newcommand{\vf}{\varphi}

\begin{document}

\begin{frontmatter}

\title{Adaptive local basis set for Kohn-Sham density
functional theory in a discontinuous Galerkin
framework I: Total energy calculation}
\author[ll]{Lin Lin} 
\ead{linlin@math.princeton.edu}

\author[jl]{Jianfeng Lu} 
\ead{jianfeng@cims.nyu.edu}

\author[ly]{Lexing Ying}
\ead{lexing@math.utexas.edu}

\author[we]{Weinan E}
\ead{weinan@math.princeton.edu}

\address[ll]{Program in Applied and Computational Mathematics, Princeton
  University, Princeton, NJ 08544.}

\address[jl]{Department of Mathematics, Courant Institute of Mathematical
  Sciences, New York University, New York, NY 10012.  }

\address[ly]{Department of Mathematics and ICES, University of Texas at
  Austin, Austin, TX 78712. }

\address[we]{Department of Mathematics and PACM, Princeton University, Princeton, NJ 08544.}

\begin{abstract}
  Kohn-Sham density functional theory is one of the most widely used
  electronic structure theories.  In the pseudopotential framework,
  uniform discretization of the Kohn-Sham Hamiltonian generally results
  in a large number of basis functions per atom in order to resolve the
  rapid oscillations of the Kohn-Sham orbitals around the nuclei.
  Previous attempts to reduce the number of basis functions per atom
  include the usage of atomic
  orbitals and similar objects, but the atomic orbitals generally
  require fine tuning in order to reach high accuracy.  We
  present a novel discretization scheme that adaptively and
  systematically builds the rapid oscillations of the Kohn-Sham
  orbitals around the nuclei as well as environmental effects into the
  basis functions.  The resulting basis functions are localized in the
  real space, and are discontinuous in the global domain. The
  continuous Kohn-Sham orbitals and the electron density are evaluated
  from the discontinuous basis functions using the discontinuous
  Galerkin (DG) framework.  Our method is implemented in parallel and
  the current implementation is able to handle systems with at least
  thousands of atoms.  Numerical examples indicate that our method
  can reach very high accuracy (less than $1$meV) with a very small number
  ($4\sim 40$) of basis functions per atom.
\end{abstract}


\end{frontmatter}


\section{Introduction}

Electronic structure theory describes the energies and distributions
of electrons, and is essential in characterizing the microscopic
structures of molecules and materials in condensed phases.  Among all
the different formalisms of electronic structure theory, Kohn-Sham
density functional theory (KSDFT)~\cite{HohenbergKohn:64, KohnSham:65}
achieves so far the best compromise between accuracy and efficiency,
and has become the most widely used electronic structure model 
for condensed matter systems.  Kohn-Sham density functional theory
gives rise to a nonlinear eigenvalue problem, which is commonly solved
using the self-consistent field iteration method~\cite{Martin:04}. In
each iteration, the Kohn-Sham Hamiltonian is constructed from a trial
electron density and is discretized into a finite dimensional matrix.
The electron density is then obtained from the low-lying
eigenfunctions, called Kohn-Sham orbitals, of the discretized
Hamiltonian.  The resulting electron density and the trial electron
density are then mixed and form a new trial electron density.  The
loop continues until self-consistency of the electron density is
reached.  An efficient algorithm therefore contains three phases:
discretization of the Hamiltonian; evaluation of the electron density
from the discretized Hamiltonian; and self-consistent iteration. In
this paper, we focus on the discretization of the Hamiltonian and the
evaluation of the electron density in the pseudopotential
framework~\cite{Martin:04}.

If space is uniformly discretized, the Kohn-Sham Hamiltonian generally
requires a basis set with a large number of degrees of freedom per
atom. For most chemical systems, the kinetic energy cutoff typically
ranges from $15$Ry to $90$Ry for standard planewave discretization in
the norm-conserving pseudopotential
framework~\cite{TroullierMartins1991}, which amounts to about $500\sim
5000$ basis functions per atom.  The required number of basis
functions per atom is even larger for uniform discretization methods
other than planewaves, such as finite difference
method~\cite{Chelikowsky:94, Alemany:04} and finite element
method~\cite{TsuchidaTsukada1995,PaskKleinFongSterne:99,
  PaskSterne2005a}.

The large number of basis functions per atom originates from the rapid
oscillation of the Kohn-Sham orbitals.  The Kohn-Sham orbitals
oscillate rapidly around the nuclei and become smooth in the
interstitial region of the nuclei.  Physical intuition suggests that
the rapid oscillations around the nuclei are inert to changes in the
environment.  A significant part of the rapid oscillations can already
be captured by the orbitals associated with isolated atoms. These
orbitals are called atomic orbitals.  Numerical methods based on
atomic orbitals or similar ideas have been designed based on this
observation~\cite{AverillEllis:73, DelleyEllis:82, Eschrig:88,
  KoepernikEschrig:99, Kenny:00, Junquera:01, Ozaki:03,
  BlumGehrkeHankeEtAl2009}. Environmental effect is not built into the
atomic orbitals directly, but can only be approximated by fine tuning
the adjustable parameters in these atomic orbitals.  The values of the
adjustable parameters therefore vary among different chemical elements
and exchange-correlation potentials, and sometimes vary among the
different ambient environment of atoms.  The quality of the atomic
orbitals are difficult to be improved systematically, but relies
heavily on the experience of the underlying chemical system.

Atomic orbitals and uniform discretization methods can be combined, as
in the mixed basis methods~\cite{Slater1937, Andersen1975, Blochl1994,
  SukumarPask:09}.  The quality of the basis functions can therefore
be systematically improved by incorporating the uniform discretization
methods. However, fine tuning the adjustable parameters is still
necessary due to the absence of the environmental effect in the basis
functions, and in certain circumstances the number of basis functions
per atom is still large.



In this paper we propose a novel discretization method to build the
environmental effects into the basis set to achieve further dimension
reduction of the basis set. The basis functions are constructed
adaptively and seamlessly from the atomic configuration in local
domains, called elements.  The basis functions are discontinuous at
the boundary of the elements, and they form the basis set used in the
discontinuous Galerkin (DG) framework.  The flexibility of the DG
framework allows us to employ these discontinuous basis functions to
approximate the continuous Kohn-Sham orbitals, and allows us to
achieve high accuracy (less than $1$meV) in the total energy calculation
with a very small number ($4\sim 40$) of basis functions per atom..
Our method is implemented in parallel with a rather general data
communication framework, and the
current implementation is able to calculate the total energy for systems
consisting of thousands of atoms. 

The discontinuous Galerkin framework has been widely used in numerical
solutions of partial differential equations (PDE) for more than four
decades, see for example 
\cite{BabuskaZlamal:73, Wheeler:78, Arnold:82,
  CockburnKarniadakisShu:00, CockburnShu:01, Arnold:02} and the
references therein. One of the main advantages of the DG method is its
flexibility in the choice of the basis functions.  The idea of
constructing basis functions adaptively from the local environment has
also been explored in other circumstances in numerical analysis such
as reduced basis method~\cite{MadayRonquist:02,
  MadayPateraTurinici:02, CancesLeBris:07, ChenHesthaven:10} and
multi-scale discontinuous Galerkin method~\cite{YuanShu:06,
  YuanShu:07, WangGuzmanShu:11} for solving PDE.
In the current context, we apply the DG algorithm to solve eigenvalue
problems with oscillatory eigenfunctions, and the basis functions are
constructed by solving auxiliary local problems numerically.

The paper is organized as follows. Section~\ref{sec:dg} introduces the
discontinuous Galerkin framework for Kohn-Sham density functional
theory. The construction of the adaptive local basis functions is
introduced in Section~\ref{sec:basis}. Section~\ref{sec:impel}
discusses implementation issues in more detail. The performance of our
method is reported in Section~\ref{sec:numer}, followed by the
discussion and conclusion in Section~\ref{sec:discussion}.


\section{Discontinuous Galerkin framework for Kohn-Sham density functional theory}
\label{sec:dg}

\subsection{Brief introduction of KSDFT}  

The Kohn-Sham energy functional in the pseudopotential
framework~\cite{Martin:04} is given by:
\begin{multline}\label{eq:KSfunc}
  E_{\KS}(\{\psi_i\}) = \frac{1}{2} \sum_{i=1}^N \int \abs{\nabla
    \psi_i}^2 \ud x + \int V_{\ext} \rho \ud x
  + \sum_{\ell} \gamma_{\ell} \sum_{i=1}^N  \abs{\int b_{\ell}^{\ast} \psi_i \ud x}^2  \\
  + \frac{1}{2} \iint \frac{\rho(x) \rho(y)}{\abs{x-y}} \ud x \ud y +
  \int\exc[\rho(x)] \ud x,
\end{multline}
where $\rho(x) = \sum_i \abs{\psi_i}^2(x)$ and the $\{\psi_i\}$'s satisfy
the orthonormal constraints:
\begin{equation}
  \int \psi_i^{\ast} \psi_j \ud x = \delta_{ij}.
\end{equation}
In \eqref{eq:KSfunc}, we have taken the Kleinman-Bylander form of
the pseudopotential~\cite{KleinmanBylander:82}. The pseudopotential
is given by
\begin{equation*}
  V_{\mathrm{PS}} = V_{\ext} + \sum_{\ell} \gamma_{\ell}
  \ket{b_{\ell}}\bra{b_{\ell}}.
\end{equation*}
For each $\ell$, $b_{\ell}$ is a function supported locally in the real
space around the position of one of the atoms , $\gamma_{\ell} = +1$ or $-1$, and we
have used the Dirac bra-ket notation. We have ignored the spin degeneracy
and have adopted the local density approximation
(LDA)~\cite{CeperleyAlder1980,PerdewZunger1981} for the
exchange-correlation functional. The proposed method can also be used
for more complicated exchange-correlation functionals and when spin
degeneracy is involved.

The Kohn-Sham equation, or the Euler-Lagrange equation associated with
\eqref{eq:KSfunc} reads 
\begin{equation}\label{eq:KSeqn}
  H_{\eff}[\rho] \psi_i = ( - \tfrac{1}{2} \Delta + V_{\eff}[\rho]
  + \sum_{\ell} \gamma_{\ell} \ket{b_{\ell}} \bra{b_{\ell}}) 
  \psi_i = E_i \psi_i, 
\end{equation}
where the effective one-body potential $V_{\eff}$ is given by 
\begin{equation}\label{eq:Veff}
  V_{\eff}[\rho](x) = V_{\ext}(x) + \int \frac{\rho(y)}{\abs{x - y}}\ud y
  + \exc'[\rho(x)].
\end{equation}
Note that \eqref{eq:KSeqn} is a nonlinear eigenvalue problem, as
$V_{\eff}$ depends on $\rho$, which is in turn determined by
$\{\psi_i\}$.  The electron density is self-consistent if both~\eqref{eq:KSeqn}
and~\eqref{eq:Veff} are satisfied. After obtaining the self-consistent
electron density, the total
energy of the system can be expressed using the eigenvalues $\{E_{i}\}$ and
$\rho$ as~\cite{Martin:04}
\begin{equation}
  E_{\text{tot}} = \sum_{i=1}^{N} E_{i} - \frac12 \iint \frac{\rho(x) \rho(y)}{\abs{x-y}} \ud x \ud y
  + \int\exc[\rho(x)] \ud x - \int \exc'[\rho(x)] \rho(x) \ud x.
  \label{eqn:Etot}
\end{equation}
The goal of Kohn-Sham density functional theory is to calculate the
total energy $E_{\text{tot}}$ and the self-consistent electron density
$\rho$ given the atomic configuration.

Numerical algorithms for Kohn-Sham density functional theory can be
broadly divided into two categories: One may try to directly minimize
the energy functional \eqref{eq:KSfunc} with respect to the Kohn-Sham
orbitals $\{\psi_i\}$ (see e.g. \cite{PayneTeterAllanEtAl1992}); one
may also try to look for a solution for \eqref{eq:KSeqn}, usually by
using the self-consistent iteration.

The self-consistent iteration goes as follows. Starting with an
initial guess $\rho_0$, one looks for a solution of \eqref{eq:KSeqn}
iteratively:
\begin{enumerate}
\item Discretization of the Hamiltonian: Determine the effective
  Hamiltonian $H_{\eff}[\rho_n]$ from the input density at the $n$-th step
  $\rho_n$;
\item Evaluation of the electron density: Obtain $\wt{\rho} = \sum_i
  \abs{\psi_i}^2$ from the effective Hamiltonian $H_{\eff}[\rho_n]$;
\item Self-consistent iteration: Determine the input density at the
  $(n+1)$-th step $\rho_{n+1}$ from $\rho_n$ and
  $\wt{\rho}$, for instance:
  \begin{equation*}
    \rho_{n+1} = \alpha \rho_n + (1 - \alpha) \wt{\rho}
  \end{equation*}
  with some parameter $\alpha$.
\item If $\norm{\rho_n - \wt{\rho}} \leq \delta$, stop; otherwise, go
  to step $(1)$ with $n \leftarrow n+1$. 
\end{enumerate}

\begin{remark} 
  The mixing step above is called linear mixing in the literature,
  which is the simplest choice. More advanced mixing
  schemes~\cite{Anderson1965,Johnson1988} can be used as well.  The
  mixing scheme used in our current implementation is the Anderson
  mixing scheme~\cite{Anderson1965}, but we will not go into the
  details of mixing schemes in this work.
\end{remark}

In this paper we focus on the discretization of the Hamiltonian and the
evaluation of the electron density. Given an effective potential $V_{\eff}$, we
find $\wt{\rho}$ from
\begin{equation}
  \wt{\rho}(x) = \sum_{i=1}^N \abs{\psi_i}^2(x),
\end{equation}
where the $\{\psi_i\}$'s are the first $N$ eigenfunctions of $H_{\eff}$.
\begin{equation}  
  H_{\eff} \psi_i = (- \tfrac{1}{2} \Delta + V_{\eff} + 
  \sum_{\ell} \gamma_{\ell} \ket{b_{\ell}} \bra{b_{\ell}}) \psi_i = E_i \psi_i.
\end{equation}
Note that the $\{\psi_i\}$'s minimize the variational problem 
\begin{equation}\label{eq:linearvar}
  E_{\eff}(\{\psi_i\}) = \frac{1}{2} \sum_{i=1}^N \int \abs{\nabla \psi_i(x)}^2
  \ud x + \int V_{\eff}(x) \rho(x) \ud x + \sum_{\ell} \gamma_{\ell} 
  \sum_{i=1}^N  \abs{ \langle b_{\ell}, \psi_i \rangle }^2,
\end{equation}
with the orthonormality constraints $\langle \psi_i, \psi_j \rangle =
\delta_{ij}$.

The evaluation of the electron density is clearly the main bottleneck in the
self-consistent iteration, which is the focus of the numerical
algorithms for Kohn-Sham density functional theory. We
consider efficient and accurate discretization for the 
evaluation of the electron density in this work.

\subsection{Discontinuous Galerkin method for KSDFT}

The discontinuous Galerkin (DG) methods have been developed for
different types of partial differential equations
\cite{BabuskaZlamal:73, Wheeler:78, Arnold:82,
  CockburnKarniadakisShu:00, CockburnShu:01, Arnold:02}.  One of the
main advantages of the DG method is its flexibility in the choice of
the approximation space, as the DG method does not require
the continuity condition of the basis functions across the interfaces of
the elements. This
flexibility is important for constructing effective discretization
schemes for Kohn-Sham
density functional theory.

We present in the following a DG method for the evaluation of the
electron density. Among the different
formalisms in the DG framework, we will use the
interior penalty method 
\cite{BabuskaZlamal:73,Arnold:82}.  The interior penalty method
naturally generalizes the variational principle~\eqref{eq:linearvar}. 

We denote by $\Omega$ the computational domain with the periodic boundary
condition, which corresponds to $\Gamma$ point sampling in the Brillouin
zone.  $\Omega$ is also referred to as the global domain in the
following discussion. Bloch boundary conditions can be taken into
account as well,
and this will appear in future publications. Let $\mc{T}$ be a collection
of quasi-uniform rectangular partitions of $\Omega$ (see
Fig.~\ref{fig:Na1Dpartition} for an example with four elements):
\begin{equation}
  \mc{T} = \{E_1, E_2, \cdots, E_M \},
\end{equation}
and $\mc{S}$ be the collection of surfaces that correspond to
$\mc{T}$.  Each $E_{k}$ is called an element of $\Omega$. For a
typical choice of partitions used in practice, the elements are chosen
to be of the same size. For example, for a crystalline material,
elements can be chosen as integer multiples of the conventional cell
of the underlying lattice. As a result, unlike the usual finite
element analysis, the element size will remain the same. \footnote{In
  the language of finite element method, we will not use the
  $h$-refinement.}

In the following discussion, we use extensively the inner products
defined as below
\begin{align}
  & \average{v, w}_E = \int_E v^{\ast}(x) w(x) \ud x, \\
  & \average{\bd{v}, \bd{w}}_S = \int_S \bd{v}^{\ast}(x)
  \cdot \bd{w}(x) \ud s(x), \\
  & \average{v, w}_{\mc{T}} = \sum_{i = 1}^M
  \average{v, w}_{E_i}, \\
  & \average{\bd{v}, \bd{w}}_{\mc{S}} = \sum_{S \in \mc{S}}
  \average{\bd{v}, \bd{w}}_S.
\end{align}
In the discontinuous Galerkin method (the interior penalty method), the
discrete energy functional corresponding to \eqref{eq:linearvar} is
given by 
\begin{multline}\label{eq:DGvar}
  E_{\DG}(\{\psi_i\}) = \frac{1}{2} \sum_{i=1}^N \average{\nabla
    \psi_i , \nabla \psi_i}_{\mc{T}} - \sum_{i=1}^N
  \average{\mean{\nabla\psi_i}, \jump{\psi_i}}_{\mc{S}}
  + \average{ V_{\eff}, \rho }_{\mc{T}} \\
  + \frac{\alpha}{h} \sum_{i=1}^N \average{\jump{\psi_i},
    \jump{\psi_i}}_{\mc{S}} + \sum_{\ell} \gamma_{\ell} \sum_{i=1}^N 
  \abs{\average{b_{\ell}, \psi_i}_{\mc{T}}}^2.
\end{multline}
Here the last term comes from the non-local terms in
Eq.~\eqref{eq:linearvar}, and $\mean{\cdot}$ and
$\jump{\cdot}$ are the average and the jump operators across surfaces, defined as
follows. For $S \in \mc{S}^{\circ}$ the set of interior surfaces, we
assume 
$S$ is shared by elements $K_1$ and $K_2$.  Denote by $n_1$ and $n_2$
the
unit normal vectors on $S$ pointing exterior to $K_1$ and $K_2$,
respectively. With $u_i = u\vert_{\partial K_i}$, $i = 1, 2$, we set
\begin{equation}
  \jump{u} = u_1 n_1 + u_2 n_2 \quad \text{on } S.
\end{equation}
For $S \in \mc{S}^{\partial}$ where $\mc{S}^{\partial}$ is the union of
the surfaces on the  boundary, we set
\begin{equation}
  \jump{u} = u n \quad \text{on } S,
\end{equation}
where $n$ is the outward unit normal. For vector-valued
function $q$, we define 
\begin{equation}
  \mean{q} = \tfrac{1}{2} (q_1 + q_2) \quad \text{on } S \in \mc{S}^{\circ},
\end{equation}
where $q_i = q \vert_{\partial K_i}$, and 
\begin{equation}
  \mean{q} = q \quad \text{on } S \in \mc{S}^{\partial}. 
\end{equation}
Note that in the current context $\mc{S} = \mc{S}^{\circ}$ since we
assume periodic boundary condition for the computational domain, and
every surface is an interior surface. The constant $\alpha$ in
\eqref{eq:DGvar} is a positive penalty parameter, which penalizes the
jumps of functions across element surfaces to guarantee stability. The
choice of $\alpha$ will be further discussed in
Section~\ref{sec:numer}.

Assume that we have chosen for each element $E_k$ a set of basis
functions $\{\varphi_{k,j}\}_{j=1}^{J_k}$, where $J_k$ is the number
of basis functions in $E_k$. We extend each $\varphi_{k,j}$ to the
whole computational domain $\Omega$ by setting it to be $0$ on the
complement set of $E_k$. Define the function space $\mc{V}$ as
\begin{equation}
  \mc{V} = \spanop\{ \varphi_{k,j},\, E_k \in \mc{T},\, j = 1, \cdots, J_k \}.
\end{equation}

We minimize \eqref{eq:DGvar} for $\{ \psi_i\} \subset \mc{V}$. The
energy functional \eqref{eq:DGvar} in the approximation space $\mc{V}$
leads to the following eigenvalue problem for $\{\psi_i\}_{i=1}^N$. For any
$v \in \mc{V}$,
\begin{multline}
  \frac{1}{2} \average{\nabla v, \nabla \psi_i}_{\mc{T}} -\frac{1}{2}
  \average{\jump{v}, \mean{\nabla \psi_i}}_{\mc{S}} -\frac{1}{2}
  \average{\mean{\nabla v}, \jump{\psi_i}}_{\mc{S}}
  +\frac{\alpha}{h} \average{\jump{v}, \jump{\psi_i}}_{\mc{S}}\\
  +\average{v, V_\eff \psi_i}_{\mc{T}} +\sum_{\ell} \gamma_{\ell}
  \average{v,b_\ell}_{\mc{T}} \average{b_\ell, \psi_i}_{\mc{T}} =
  \lambda_i \average{v, \psi_i}_{\mc{T}}.
\end{multline}
Setting $v = \varphi_{k', j'}$ and
\begin{equation}
  \psi_i = \sum_{E_k \in \mc{T}} \sum_{j=1}^{J_k} c_{i; k, j}
  \varphi_{k, j},
\end{equation}
we arrive at the following linear system
\begin{multline}
  \sum_{k,j} \biggl( \frac{1}{2}\average{\nabla \vf_{k',j'}, \nabla
    \vf_{k,j}}_{\mc{T}} -\frac{1}{2}\average{\jump{\vf_{k',j'}},
    \mean{\nabla \vf_{k,j}}}_{\mc{S}} \\
  -\frac{1}{2}\average{\mean{\nabla \vf_{k',j'}},
    \jump{\vf_{k,j}}}_{\mc{S}} +\frac{\alpha}{h}
  \average{\jump{\vf_{k',j'}}, \jump{\vf_{k,j}}}_{\mc{S}}
  +\average{\vf_{k',j'}, V_\eff \vf_{k,j}}_{\mc{T}} \\
  + \sum_\ell \gamma_{\ell} \average{\vf_{k',j'},b_\ell}_{\mc{T}}
  \average{b_\ell, \vf_{k,j}}_{\mc{T}} \biggr) c_{i;k,j} =
  \lambda_{i} \sum_{k,j}
  \average{\vf_{k',j'}, \vf_{k,j}} c_{i;k,j}.
\label{eq:DGimpl}
\end{multline}
We define $A$ to be the matrix with entries given by the expression in
the parentheses, $B$ to be the matrix with entries
$\average{\vf_{k',j'}, \vf_{k,j}}$, and $c_i$ to be the vector with
components $(c_{i;k,j})_{k,j}$, we have the following simple form of
generalized eigenvalue problem
\[
A c_i = \lambda_i B c_i
\]
for $i=1,2,\ldots, N$.  
Following the standard terminologies in the finite element method, we
call $A$ the (DG) stiffness matrix, and $B$ the (DG) mass matrix.  In
the special case when the DG mass matrix $B$ is
equal to the identity matrix, we have a standard eigenvalue problem 
$A
c_i = \lambda_i c_i$.  Once $\{c_i\}$ are available, the electron
density is calculated by
\begin{equation}\label{eq:newdensity}
  \wt{\rho} = \sum_{i=1}^N \sum_{E_k\in\mc{T}} \abs{\sum_{j=1}^{J_k} 
  c_{i; k, j} \varphi_{k, j}}^2.
\end{equation}

\section{Basis functions adapted to the local environment}\label{sec:basis}

The proposed framework in the last section is valid for any choice of
basis functions. To improve the efficiency of the algorithm, it is
desirable to use less number of basis functions while maintaining the
same accuracy. In order to achieve this goal, the choice of the
functions $\{\varphi_{k,j}\}$ is important. In this section, we
discuss a way to construct the basis functions $\{\varphi_{k,j}\}$
that are adapted to the local environment.

\begin{figure}[ht]
  \begin{center}
    \subfloat[$\rho$]{\includegraphics[width=0.40\textwidth]{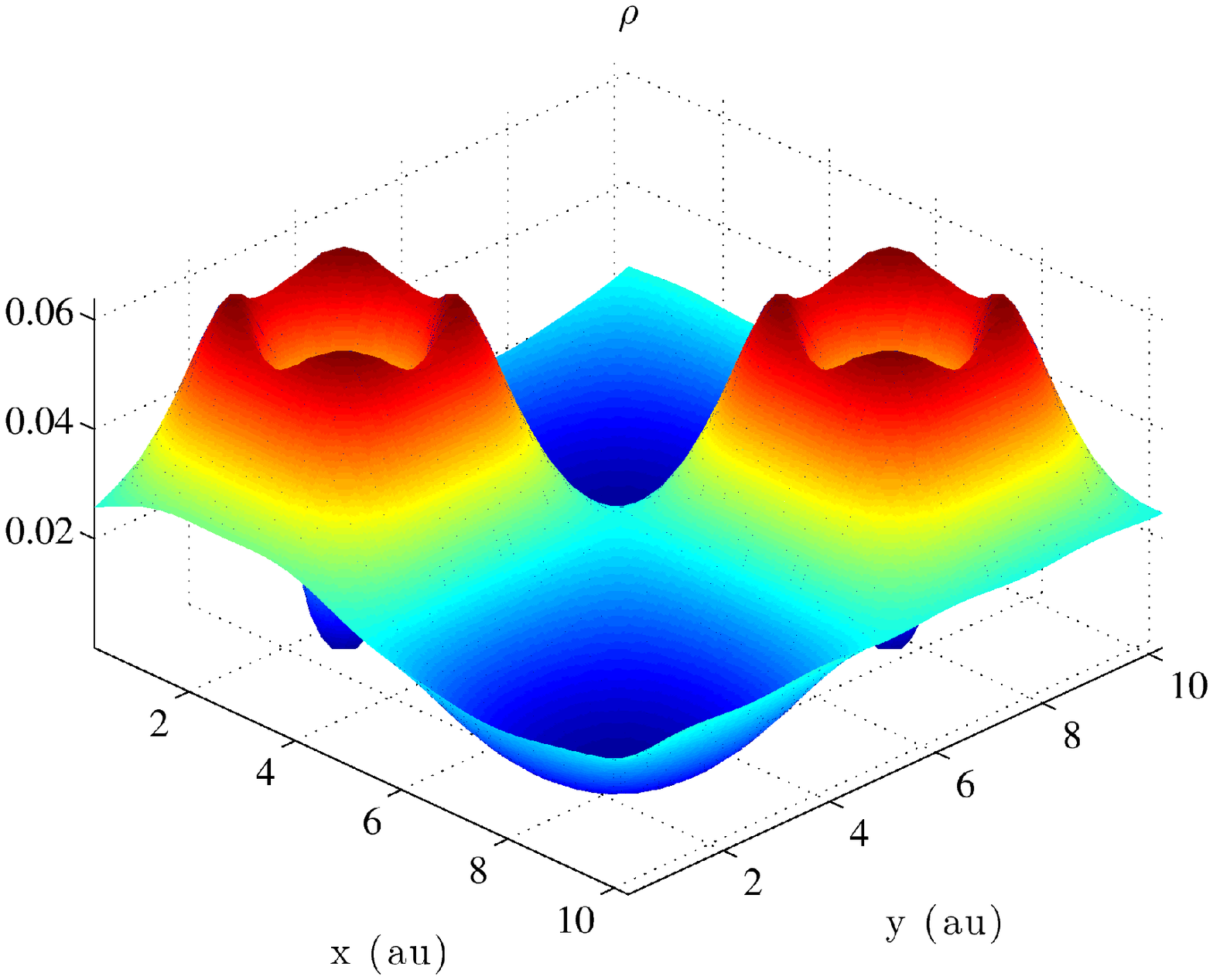}}
    \subfloat[$\norm{\nabla\rho}_2$]{\includegraphics[width=0.40\textwidth]{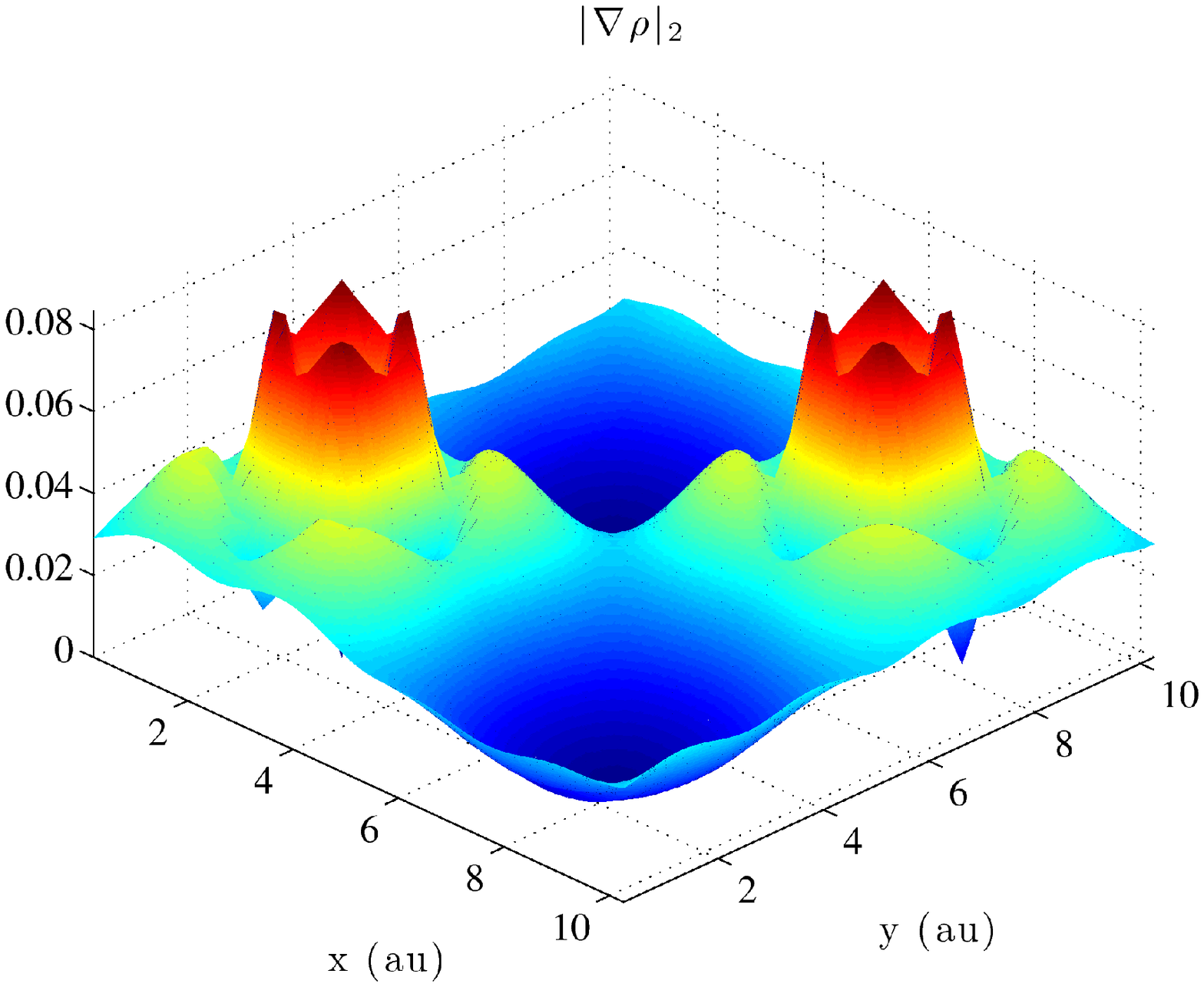}}
  \end{center}
  \caption{(color online) The electron density (a) and the norm of the
  gradient of the electron density (b) on a $(001)$ slice of
  a mono-crystalline silicon system passing through two Si atoms. The
  two Si atoms are located at $(2.57,2.57)$~au and at
  $(7.70,7.70)$~au in this plane, respectively. Even in the
  pseudopotential framework, the electron density
  is much more oscillatory around the nuclei of the Si atoms and is smooth in the
  interstitial region. } 
  \label{fig:sing} 
\end{figure}

The starting point is the following observation. The Kohn-Sham
orbitals $\{\psi_i\}$ exhibit oscillatory behavior around the nuclei.  In a
full electron calculation, the nuclei charge density is the summation
of delta functions located at the positions of the nuclei (or
numerical delta function after discretization) and the Kohn-Sham
orbitals have cusp points at the positions of the atoms. In the
pseudopotential framework which involves only valence electrons, one
can still see that the Kohn-Sham orbitals and the electron density are
much more oscillatory near the atom cores than in the interstitial
region, as illustrated in Fig.~\ref{fig:sing}.  In the setting of
real space method or planewave method, in order to resolve the
Kohn-Sham orbitals around the atom cores where the derivatives of
Kohn-Sham orbitals become large one has to use a uniform fine
mesh. Therefore, the number of mesh points becomes huge even for a
small system. This makes the electronic structure calculation
expensive.

In order to reduce the cost, we note that the Kohn-Sham orbitals are
smooth away from the atoms and the uniform fine discretization is not
efficient enough. Adaptive refinement
techniques can be used to improve the efficiency by reducing the number
of basis functions per atoms.  Techniques of this type include finite
element based adaptive mesh refinement
method~\cite{FattebertHornungWissink2007}, finite volume based adaptive
mesh refinement method, and multiresolution basis
functions~\cite{Arias1999,GenoveseNeelovGoedeckerEtAl2008,HarrisonFannYanaiEtAl2004},
to name a few.  Our approach builds the oscillatory behavior
the Kohn-Sham orbitals near the atom cores into the basis functions.
Hence, a small number of basis functions are enough to characterize the
Kohn-Sham orbitals. This idea is not entirely new. For
example, the philosophy of pseudopotential techniques is quite
similar, though the reduction is done at the analytic level. On the
side of numerical methods, the ideas behind atomic orbital basis or
numerical atomic basis are closely
related~\cite{Eschrig:88,BlumGehrkeHankeEtAl2009}.

The main difference from the previous approaches is that instead of
predetermining basis functions based on the information from isolated
atoms, our approach builds the information from the local environment
into the basis functions as well. Thanks to the flexibility of the
discontinuous Galerkin framework, this can be done in a seamless and
systematic way.  The basis functions form a complete basis set in the
global domain $\Omega$. The basis set is therefore efficient, and at
the same time the accuracy can be improved systematically. This is an
important difference between this approach and the previous methods
along the same line.

The basis functions $\{\varphi_{k, j}\}$ are determined as follows.
Given the partition $\mc{T}$ and the effective potential $V_{\eff}$, let
us focus on the construction of $\{\varphi_{k, j}\}$, $j = 1, \cdots,
J_k$ for one
element $E_k \in \mc{T}$. As discussed above, our approach is to adapt
$\{\varphi_{k,j}\}$ to the local environment in $E_k$.

For each element $E_k$, we take a region $Q_k \supset E_k$. $Q_k$ is
called the extended element associated with the element $E_k$.  The
set $Q_{k}\backslash E_{k}$ is called the buffer area.  We will choose
$Q_{k}$ which extends symmetrically along the $\pm x(y,z)$ directions
from the boundary of $E_k$.  The length of the buffer area extended
beyond the boundary of $E_{k}$ along the $\pm x(y,z)$ direction is
called the ``buffer size along the $x(y,z)$ direction''.  We restrict
the effective Hamiltonian on $Q_k$ by assuming the periodic boundary
condition on $\partial Q_k$ and denote by $H_{\eff, Q_k}$ the
restricted Hamiltonian.  $H_{\eff, Q_k}$ is discretized and
diagonalized, and the corresponding eigenfunctions are denoted by
$\{\widetilde{\varphi}_{k, j}\}$, indexed in increasing order of the
associated eigenvalues. We restrict the first $J_k$ eigenfunctions
$\{\widetilde{\varphi}_{k, j}\}$ from $Q_k$ to $E_k$, denoted by
$\{\varphi_{k, j}\}$.  Each $\varphi_{k, j}$ is therefore defined
locally on $E_k$.  As discussed before we extend each $\varphi_{k,j}$
to the global domain $\Omega$ by setting the value to be $0$ on the
complement of $E_k$. The resulting functions, still denoted by
$\{\varphi_{k,j}\}$ are called the adaptive local basis functions.
Numerical result suggests that we can take very small $J_k$ to achieve
chemical accuracy.


The reason why we choose the periodic boundary condition on $Q_k$ for the
restriction $H_{\eff, Q_k}$ is twofold. On one hand, the periodic
boundary condition captures better the bulk behavior of the system
(than the Dirichlet boundary condition for example); On the other hand,
the periodic boundary condition makes the solution of $H_{\eff, Q_k}$
more easily adapted to existing DFT algorithms and packages, as most of
them can treat periodic boundary conditions. Other choices such as
the Neumann boundary condition are possible, and the optimal choice of
boundary conditions remains to be an open question.

The basis functions constructed from the buffer region capture well
the local singular behavior of Kohn-Sham orbitals near the
nuclei. Hence, the approximation space formed by $\{\varphi_{k,j}\}$
gives an efficient and accurate discretization to the problem, as will
be illustrated by numerical examples in Section~\ref{sec:numer}.  Note
that the $\{\widetilde{\varphi}_{k,j}\}$'s are the eigenfunctions of
the self-adjoint operator $H_{\eff, Q_k}$ on $Q_k$, and therefore form
a complete basis set on $Q_k$ when $J_k \to \infty$. This implies that
after restriction, the functions $\{\varphi_{k,j}\}$ also form a
complete basis set on $E_k$ as $J_k \to \infty$.  The accuracy can
therefore be systematically improved in the electronic structure
calculation.





Eq.~\eqref{eq:DGimpl} proposes a generalized eigenvalue problem. From
numerical point of view it would be more efficient if we can choose
$\{\varphi_{k,j}\}$ such that the DG mass matrix is an identity matrix
and that Eq.~\eqref{eq:DGimpl} becomes a standard eigenvalue problem.
Moreover, as $J_{k}$ increases, the basis functions
$\{\varphi_{k,j}\}$ can become degenerate or nearly degenerate, which
increases the condition number. Both problems can be solved at the
same time by applying a singular value decomposition (SVD) filtering
step, resulting in an orthonormal basis set $\{\varphi_{k,j}\}$:
\begin{enumerate}
\item For each $k$, form a matrix $M_k = ( \varphi_{k, 1}, \varphi_{k, 2},
  \cdots, \varphi_{k, J_k} ) $ with $\varphi_{k,j}$;
\item Calculate SVD decomposition $U D V^{\ast} = M_k$,
  \[D = \diag(\lambda_{k, 1}, \lambda_{k, 2}, \cdots, \lambda_{k,
    J_k}), \] where $\lambda_{k, j}$ are singular values of $M_k$
  ordered decreasingly in magnitude;
\item For a threshold $\delta$, find $\wt{J}_k$ such that
  $\abs{\lambda_{k, \wt{J}_k}} > \delta$ and $\abs{\lambda_{k,
      \wt{J}_k+1}} < \delta$ ($\wt{J}_k = J_k$ if all singular values
  are larger than the threshold). Take $U_{j}$ be the
  $j$-th column of $U$, $j = 1, \cdots, \wt{J}_k$;
\item Set $J_k \leftarrow \wt{J}_k$ and $\varphi_{k, j} \leftarrow
  U_{k,j}$ for $j = 1, \cdots, \wt{J}_k$.
\end{enumerate}
\begin{remark}
  Although the threshold $\delta$ can avoid numerical degeneracy of
  the basis functions, the numerical degeneracy is not observed for
  the cases studied in section~\ref{sec:numer}.  In other words, we
  will take $\delta = 0$, $J_k = \wt{J}_k$.
\end{remark}

After constructing the basis functions $\{\varphi_{k, j} \}$, we then
apply the discontinuous Galerkin framework to solve $\{\psi_{i}\}$ and
hence $\rho$ corresponding to $H_{\eff}$. We summarize the overall
algorithm as follows:

\begin{enumerate}
  \item Set $n = 0$, let $\mc{T}$ be a partition of $\Omega$ into
    elements, and $\rho_{0}$ be an initial trial electron density;
\item Form the effective potential $V_{\eff}[\rho_n]$ and the
  effective Hamiltonian $H_{\eff}[\rho_n]$;
\item For each element $E_k \in \mc{T}$, calculate the eigenfunctions
  corresponding to the Hamiltonian $H_{\eff,Q_k}$ on the extended
  element $Q_k$, and obtain the orthonormal adaptive local basis
  functions $\{\varphi_{k,j}\}$;
\item Solve \eqref{eq:DGimpl} to obtain the coefficients
  $\{c_{i;k,j}\}$ for the Kohn-Sham orbitals and reconstruct the
  electron density $\wt{\rho}$ by \eqref{eq:newdensity};
\item Mixing step: Determine $\rho_{n+1}$ from $\rho_n$ and
  $\wt{\rho}$. If $\norm{\rho_n - \wt{\rho}} \leq \delta$, stop; otherwise, go
  to step $(2)$ with $n \leftarrow n+1$. 
\end{enumerate}

We remark that due to the flexibility of the DG framework one can
supplement the functions $\{\varphi_{k,j}\}$ constructed above by
other functions in $E_k$, such as local polynomials in $E_k$, Gaussian
functions restricted to $E_k$, and other effective basis functions
based on physical and chemical intuition. From practical point of
view, we find that the adaptive basis set constructed above already
achieves satisfactory performance.


\section{Implementation details}\label{sec:impel}

This section explains the implementation details for the above
algorithm. Specialists of the DG methods can skip this section and go
directly to the numerical results in Section~\ref{sec:numer}. This
section is mostly written for the readers who are less familiar with the
DG implementation.

\subsection{Grids and Interpolation}

The above algorithm involves three types of domains: the global
domain $\Omega$, the extended elements $\{Q_k\}$, and the elements
$\{E_k\}$. Quantities defined on these domains are discretized with
different types of grids.
\begin{itemize}
\item On $\Omega$, the quantities such as $\rho$ and $V_\eff$ are
  discretized with a uniform Cartesian grid with a spacing fine enough
  to capture the singularities and oscillations in these quantities.
\item The grid on $Q_k$ is simply the restriction of the uniform grid
  of $\Omega$ on $Q_k$. This is due to the consideration that all
  quantities on $Q_k$ are treated as periodic and hence a uniform grid
  is the natural choice.
\item The grid on $E_k$ is a three-dimensional Cartesian
  Legendre-Gauss-Lobatto (LGL) grid in order to accurately carry out
  the operations of the basis functions $\{\varphi_{k,j}\}$ such as
  numerical integration and trace operator for each element $E_k$.
\end{itemize}
Transferring various quantities between these three grids requires the
following interpolation operators.
\begin{itemize}
\item $\Omega$ to $Q_k$. This is used when we restrict the density
  $\rho_n$ and the effective potential $V_\eff$ to the extended element $Q_k$.
  Since the grid on $Q_k$ is the restriction of the grid on $\Omega$,
  this interpolation operator simply copies the required values.
\item $Q_k$ to $E_k$. This is used when one restricts
  $\{\widetilde{\varphi}_{k,j}\}$ and their derivatives to $E_k$. As
  the grid on $Q_k$ is uniform, the interpolation is done by Fourier
  transform. Due to the fact that both grids are Cartesian, the
  interpolation can be carried out dimension by dimension, which
  greatly improves the efficiency.
\item $E_k$ to $\Omega$. This is used when one assembles the Kohn-Sham
  orbitals $\{\psi_i\}$ from the coefficients $\{c_{i;k,j}\}$ of the
  elements. The interpolation from the LGL grid to the uniform grid is
  done by Lagrange interpolation, again carried out dimension by
  dimension. Averaging is performed for the grid points of $\Omega$
  shared by multiple elements.
\end{itemize}

The non-local pseudopotentials are used both in solving
$\{\wt{\varphi}_{k,j}\}$ on each $Q_k$ and in the numerical
integration step on the LGL grid of each $E_{k}$. In our
implementation, the non-local pseudopotentials are directly generated
in real space on $Q_{k}$ and on $E_{k}$ without further interpolation
between the grids.

\subsection{Implementation of the discontinuous Galerkin method}

We use planewaves in each extended element $Q_k$ to discretize the
local effective Hamiltonian $H_{\eff, Q_k}$ and the LOBPCG
algorithm~\cite{Knyazev:01} with the preconditioner proposed
in~\cite{TeterPayneAllan1989} to diagonalize the discretized
Hamiltonian. The resulting eigenfunctions
$\{\wt{\varphi}_{k,j}\}_{j=1}^{J_k}$ of $H_{\eff, Q_k}$ are restricted
to $E_k$ and interpolated onto its LGL grid.
Within the SVD filtering step, the
inner product that we adopt is the discrete weighted $\ell_2$ product
with the LGL weights inside $E_k$. The main advantage of the SVD
filtering step is that the discontinuous Galerkin method results in 
a standard eigenvalue problem.

The assembly of the DG stiffness matrix follows \eqref{eq:DGimpl} and
consists of the following steps.
\begin{itemize}
\item For the first term $\frac{1}{2}\average{\nabla \vf_{k',j'},
    \nabla \vf_{k,j}}_{\mc{T}}$ and the fifth term
  $\average{\vf_{k',j'}, V_\eff \vf_{k,j}}_{\mc{T}}$, their
  contributions are non-zero only when $k=k'$ since otherwise two
  basis functions have disjoint support. Hence, for each fixed $k$, we
  compute $\average{\nabla \vf_{k,j'}, \nabla \vf_{k,j}}_{E_k}$ and
  $\average{\vf_{k,j'}, V_\eff \vf_{k,j}}_{E_k}$. The integration is
  done numerically using the LGL grid on $E_k$. Part of the stiffness
  matrix corresponding to these two terms clearly has a block diagonal
  form.
\item For the second, third, and fourth terms of \eqref{eq:DGimpl},
  one needs to restrict basis functions and their derivatives to
  element faces. As the one-dimensional LGL grid contains the
  endpoints of its defining interval, this is done simply by
  restricting the values of the three-dimensional LGL grid to the
  element faces. One then calculates these three terms using numerical
  integration on these resulting two-dimensional LGL grids. Since the
  integral is non-zero only when $E_k$ and $E_{k'}$ are the same
  element or share a common face, part of the stiffness matrix
  corresponding to these three terms is again sparse.
\item The last term of \eqref{eq:DGimpl} is $\sum_\ell \gamma_{\ell}
  \average{\vf_{k',j'},b_\ell}_{\mc{T}} \average{b_\ell,
    \vf_{k,j}}_{\mc{T}}$. The integration is again approximated using
  the LGL grids of the elements. Notice that the contribution is
  non-zero iff $\vf_{k',j'}$ and $\vf_{k,j}$ overlap with the support
  of a common $b_\ell$. Since each $b_\ell$ is localized around a
  fixed atom, $\vf_{k,j}$ and $\vf_{k',j'}$ need to be sufficiently
  close for this term to be non-zero. As a result, part of the
  stiffness matrix corresponding to this last term is also sparse.
\end{itemize}
Though the DG stiffness matrix $A$ is sparse, this property is not yet
exploited in the current implementation. The eigenvalues and
eigenvectors of the DG stiffness matrix are calculated using the
\textsf{pdsyevd} routine of ScaLAPACK by treating it as a dense
matrix. We plan to replace it with more sophisticated solvers that
leverage the sparsity of $A$ in future.

\subsection{Parallelization}

Our algorithm is implemented fully in parallel for message-passing
environment. To simplify the discussion, we assume that the number of
processors is equal to the number of elements. It is then convenient
to index the processors $\{P_k\}$ with the same index $k$ used for the
elements. In the more general setting where the number of elements is
larger than the number of processors, each processor takes a couple of
elements and the following discussion will apply with only minor
modification. Each processor $P_k$ locally stores the basis functions
$\{\vf_{k,j}\}$ for $j=1,2,\ldots,J_k$ and the unknowns
$\{c_{i;k,j}\}$ for $i=1,2,\ldots,N$ and $j=1,2,\ldots,J_k$. We
further partition the non-local pseudopotentials $\{b_\ell\}$ by
assigning $b_\ell$ to the processor $P_k$ if and only if the atom
associated to $b_\ell$ is located in the element $E_k$.


The eigenfunctions of the local Hamiltonian $H_{\eff,Q_k}$ are
calculated on each processor $P_k$. In order to build the local
Hamiltonian $H_{\eff,Q_k}$, the processor $P_k$ needs to access all the
non-local pseudopotentials of which the associated atoms are located in
$Q_k$.  This can be achieved by communication among $E_k$ and its nearby
elements. Once these pseudopotentials are available locally, the
eigenfunctions of $H_{\eff,Q_k}$ are computed in parallel without any
extra communication between the processors.

The parallel implementation of the DG solve is more
complicated:
\begin{itemize}
\item For the calculation of the first and the fifth terms of the
  DG stiffness matrix $A$ in Eq.~\eqref{eq:DGimpl}, each processor $P_k$
  performs numerical integration on $E_k$. Since the local basis
  functions $\{\vf_{k,j}\}$ are only non-zero on $E_k$, this step is
  carried out fully in parallel.
\item To calculate the second, third, and fourth terms, each processor
  $P_k$ computes the surface integrals restricted to the {\em left,
    front,} and {\em bottom} faces of $E_k$. This requires the basis
  functions of the left, front, and bottom neighboring elements.
\item To calculate the sixth term, each processor $P_k$ computes the
  parts associated with the non-local pseudopotentials $\{b_\ell\}$
  located on $P_k$. This requires the access to the basis functions of
  all elements that overlap with $b_\ell$.
\end{itemize}
To summarize, each processor $P_k$ needs to access the basis functions
from its neighboring elements and from the elements that overlap with
the support set of the non-local pseudopotentials located on the
elements associated with $P_k$. Due to the locality of the non-local
pseudopotentials, these elements are geometrically close to $P_k$. Since
the size of the elements is generally equal to or larger than one unit
cell, the support set of the non-local pseudopotentials are generally
within the range of the neighboring elements. Therefore, the number of
the non-local basis functions required by $P_k$ is bounded by a small
constant times the typical number of the basis functions in an element.

The use of the \textsf{pdsyevd} routine of ScaLAPACK for solving the
eigenvalue problem ~\eqref{eq:DGimpl} results in another source of
communication. ScaLAPACK requires $A$ to be stored in its block cyclic
form and this form is quite different from the distribution in which
the DG stiffness matrix is assembled (as mentioned above). As a
result, one needs to redistribute $A$ into this block cyclic form
before calling \textsf{pdsyevd} and then redistribute the
eigenfunctions afterwards.

In order to support these two sources of data communication, we have
implemented a rather general communication framework that only
requires the programmer to specify the desired non-local data. This
framework then automatically fetches the data from the processors that
store them locally. The actual communication is mostly done using
asynchronous communication routines \textsf{MPI\_Isend} and
\textsf{MPI\_Irecv}.



\section{Numerical examples}\label{sec:numer}

In order to illustrate how our method works in practice, we present
numerical results for the ground state electronic structure
calculation, using sodium (Na) and silicon (Si) as representative
examples for metallic and insulating systems, respectively.  We find
that very high accuracy (less than $10^{-6}$ au per atom) can be
achieved by using only a small number of adaptive local basis
functions.  Because of the small number of basis functions per atom,
the DG scheme already exhibits significant speedup in computational
time for a small system containing $128$ Na atoms.  We demonstrate
that the current implementation is able to solve systems with
thousands of atoms, and that the algorithm has a potential to be
applied to much larger systems with a more advanced implementation.

This section is organized as follows: section~\ref{subsec:setup}
introduces the setup of the test systems and how the error is quantified.  
Section~\ref{subsec:1Dperiodic} applies
the adaptive local basis functions to disordered quasi-1D sodium and
silicon system, followed by the result for the disordered quasi-2D
and bulk 3D systems in section~\ref{subsec:2D3D}.
We discuss the effect of the penalty parameter $\alpha$ in
section~\ref{subsec:penalty}.  Finally we demonstrate the computational
performance of our parallel implementation of the adaptive local basis
functions in section~\ref{subsec:efficiency}.

\subsection{Setup}\label{subsec:setup}




We use the local density approximation
(LDA)~\cite{CeperleyAlder1980,PerdewZunger1981} for the
exchange-correlation functional, and Hartwigsen-Goedecker-Hutter (HGH)
pseudopotential~\cite{HartwigsenGoedeckerHutter1998} with the local and
non-local pseudopotential fully implemented in the real
space~\cite{PaskSterne2005}.  All quantities are reported in atomic
units (au).  
All calculations are carried out on the Hopper system
maintained at National Energy Research Scientific Computing Center
(NERSC). Each compute node on Hopper has $24$ processors (cores) with
$32$ gigabyte (GB) of memory (1.33 GB per core). 

\begin{figure}[ht]
  \begin{center}
    \includegraphics[width=0.7\textwidth]{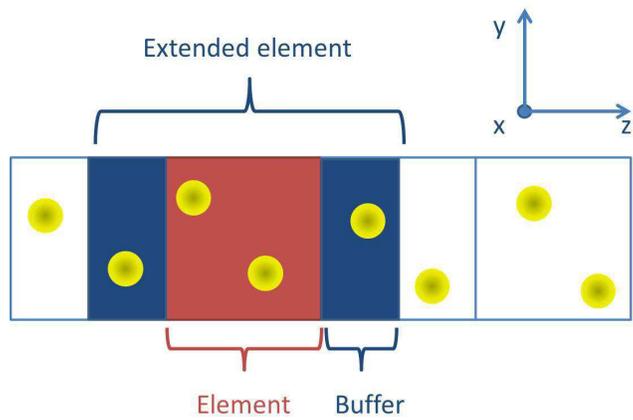}
  \end{center}
  \caption{(color online) A quasi-1D disordered Na system with 8 atoms
  extended along the $z$ direction, viewed along the $x$ direction.  The length of each empty box is
  equal to the lattice constant for the perfect Na crystal.  The red
  area represents one of the elements $E_{2}$. The corresponding
  extended element $Q_2$ consists of both the red area and the blue area
  (buffer). The buffer size is $0.50$ (in the unit of lattice constant)
  along the $z$ direction, and is $0.0$ along the $x$ and $y$
  directions.  } 
  \label{fig:Na1Dpartition}
\end{figure}

The performance of the adaptive local basis functions are tested using
Na and Si as representative examples for simple metallic and
insulating systems, respectively.  The crystalline Na has a body
centered cubic (bcc) unit cell, with 2 atoms per cell and a lattice
constant of $7.994$~au. The crystalline Si has a diamond cubic unit
cell, with 8 atoms per cell and a lattice constant of
$10.261$~au. Each atomic configuration in the following tests is
obtained by forming a supercell consisting $m\times n\times p$ unit
cells with perfect crystal structure, and a random displacement
uniformly distributed in $[-0.2, 0.2]$ au is then applied to each
Cartesian coordinate of each atom in the supercell.  The resulting
atomic configuration is therefore mildly disordered in order to avoid
the possible cancellation of errors for the case of perfect
crystalline systems.  A system is called quasi-1D if $1=m=n<p$,
quasi-2D if $1=m<n=p$, and 3D bulk if $1<m=n=p$, respectively.  In all
the tests below, the element is chosen to be the (conventional) unit
cell of the lattice.  Fig.~\ref{fig:Na1Dpartition} shows how a
quasi-1D Na system with 8 atoms extended along the $z$ direction are
partitioned in order to generate adaptive local basis functions.  The
global domain is partitioned into 4 elements $\{E_{k}\}_{k=1}^{4}$
with $2$ atoms per element.  The red area represents one of the
elements $E_{2}$, and the corresponding extended element $Q_{2}$
consists of both the red area and the blue area (buffer).  We recall
that the buffer size along the $x(y,z)$ direction refers to the length
of the buffer area extended beyond the boundary of the element $E_{k}$
along the $x(y,z)$ direction. The unit of buffer size is the lattice
constant for the perfect crystalline system.
Fig.~\ref{fig:Na1Dpartition} shows the case with the buffer size of
$0.50$ along the $z$ direction, and $0.0$ along the $x$ and $y$
directions.

We quantify the error of the adaptive local basis functions the error
of the total energy per atom which is defined as follows.  First, the
electronic structure problem is solved by using planewaves on the
global domain starting from a random initial guess of the electron
wavefunctions.  The total energy after reaching self-consistency is
denoted by $E_{GLB}$.  Then, the same electronic structure problem is
solved by the DG formulation starting from a random initial guess of
the adaptive local basis functions on each element.  The total energy
after reaching self-consistency is denoted by $E_{DG}$.  The global
domain calculation and the DG calculation using adaptive local basis
functions are therefore completely independent, and the error of the
total energy per atom is defined to be
$\abs{E_{GLB}-E_{DG}}/N_{atom}$.  For simplicity only $\Gamma$ point
is used in the Brillouin zone sampling. The proposed method can be
easily generalized to $k$-point sampling.  $10$ LOBPCG iterations are
used in each SCF iteration for the global domain calculation, and $3$
LOBPCG iterations are used in each SCF iteration for generating the
adaptive local basis functions in the DG calculation.  A small number
of LOBPCG iterations is already sufficient, since the electron
wavefunctions in the global domain calculation and the adaptive local
basis functions in the DG calculations at the end of each SCF
iteration can be reused as the initial guess in the consequent SCF
iteration for continuous refinement.  Anderson mixing is used for the
SCF iteration with a fictitious electron temperature set to be $2000$
K to facilitate the convergence of the SCF iteration.

The grid spacing for the global domain calculation is
$0.4$ au for Na and $0.32$ au for Si.  This translates to a grid of size
$20\times 20\times 20$ to discretize one Na unit cell and a grid of size
$32\times 32\times 32$ to discretize one Si unit cell.  The
Legendre-Gauss-Lobatto (LGL) grid for each element is $20\times
20\times 20$ for Na and $40\times 40\times 40$ for Si.  The LGL grid is
only used for the purpose of numerical integration in the assembly
process of the DG matrix. We remark that this grid is denser than what
is commonly used for the electronic structure calculation for three
reasons: 1) The HGH pseudopotential used in the present calculation is
more stiff than many other pseudopotentials such as the Troullier-Martins
pseudopotential~\cite{TroullierMartins1991}; 2) The potentials and
wavefunctions are represented in the real space rather than in the Fourier
space; 3) Most importantly, a dense grid in the real space is needed in both
global domain calculations and DG calculations in order to reliably reflect the
error of the total energy per atom.  

We remarked in the end of section~\ref{sec:basis} that the DG framework
is very flexible and can incorporate not only the adaptive local basis
functions but also other basis functions such as local polynomials.
In practice we find that the adaptive local basis functions are
computationally more efficient than polynomials.  Therefore in the
following discussion only adaptive local functions will be used in the
basis set.  The number of adaptive local functions per atom is also
referred to as the degrees of freedom (DOF) per atom.

\subsection{Disordered Quasi-1D System}\label{subsec:1Dperiodic}

\begin{figure}[ht]
  \begin{center}
    \subfloat[Quasi-1D Na]{\includegraphics[width=0.45\textwidth]{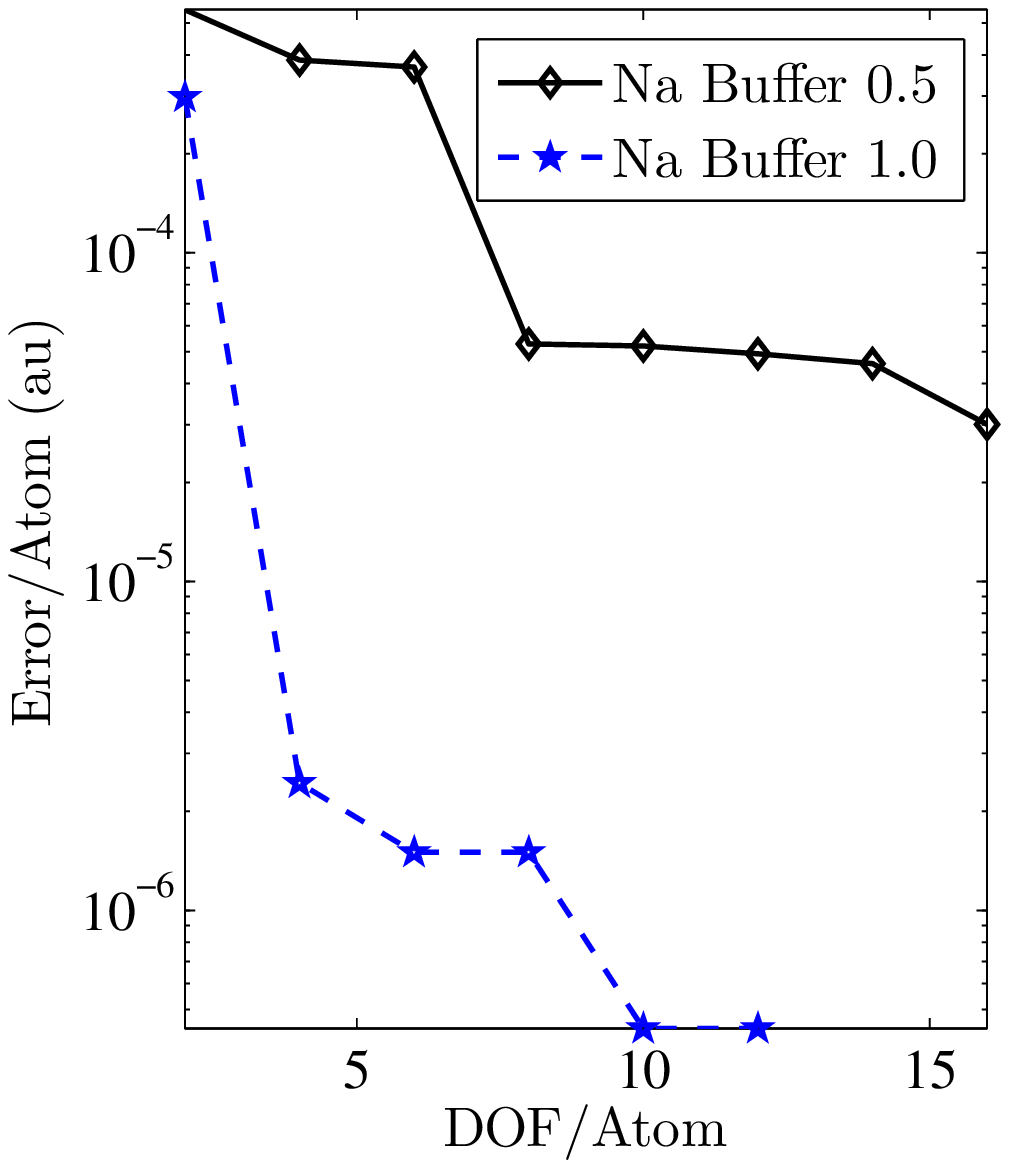}}
    \subfloat[Quasi-1D Si]{\includegraphics[width=0.45\textwidth]{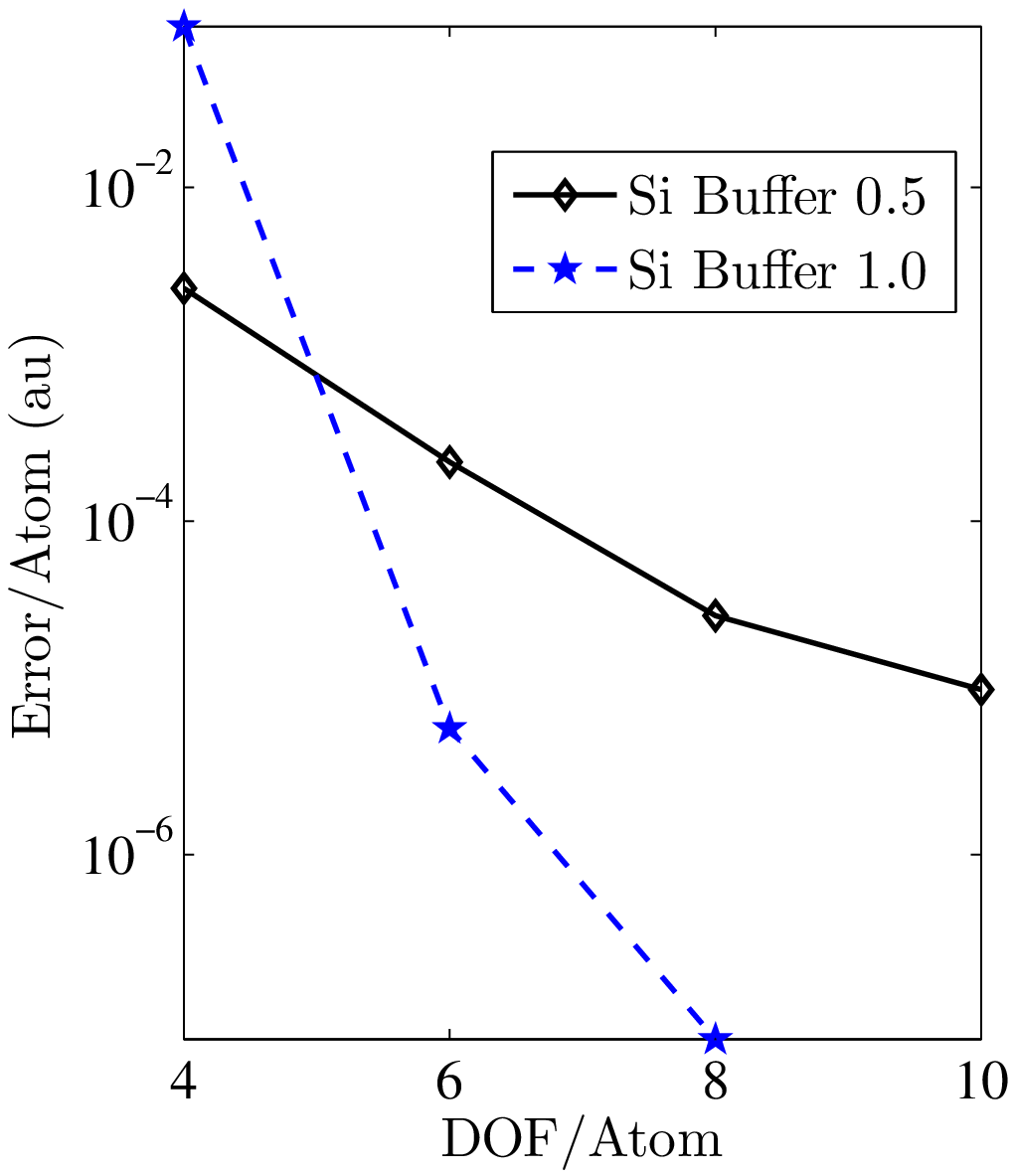}}
  \end{center}
  \caption{(color online) (a) The error of the total energy per atom 
  (the $y$ axis, plotted in log-scale) for a disordered quasi-1D sodium system consisting of
  $8$ atoms, with
  respect to the number of adaptive local basis functions per atom (the
  $x$ axis). The buffer sizes are chosen to be $0.50$ (black diamond
  with solid line), and  $1.00$ (blue star with dashed line).  (b) The
  error of the total energy per atom (the $y$axis, plotted in log-scale) for
  a disordered quasi-1D silicon system consisting of $32$ atoms, with
  respect to the number of adaptive local basis functions per atom (the
  $x$ axis). The legend is the same as in (a).  }
  \label{fig:NaSienrich}
\end{figure}

The error of the total energy per atom with respect to different
buffer sizes and different numbers of basis functions per atom (DOF
per atom) is illustrated for the disordered quasi-1D sodium system in
Fig.~\ref{fig:NaSienrich} (a) and for the disordered quasi-1D silicon
system in Fig.~\ref{fig:NaSienrich} (b).  The penalty parameter
$\alpha$ is $20$.  In both cases, the error decreases systematically when
the buffer size and the number of basis functions per atom increase.
For Na, the error of the total energy per atom is already below
$10^{-3}$ au using as few as $4$ basis functions per atom with a small
buffer of size $0.50$ (black diamond with solid line).  When the
buffer size is increased to $1.00$ (blue star with dashed line), the
error of the total energy per atom is $4.3\times 10^{-7}$ au or $0.01$
meV using $10$ adaptive local basis functions per atom.

Similar behavior is found for the silicon system.  For a small buffer of
size $0.50$ (black diamond with solid line), the error of the total
energy per atom
is $2.3\times 10^{-4}$ au with $6$ basis functions per atom.  For the
buffer of size $1.00$ (blue star with dashed line), the error of the total energy per atom is
$7.8\times 10^{-8}$ au or $0.002$ meV using as few as $8$ basis
functions per atom.  
Physical intuition suggests that the minimum number of basis functions
is $4$, which reflects one $2s$ and three $2p$ atomic orbitals.
$20\sim 40$ number of basis functions per atom is generally required
to achieve good accuracy if Gaussian type orbitals or numerical atomic
orbitals are to be used~\cite{BlumGehrkeHankeEtAl2009}.  Therefore for
the quasi-1D systems tested here, our algorithm achieves nearly the optimal
performance in terms of the number of basis functions per atom.

\begin{figure}[ht]
  \begin{center}
    \subfloat[Na]{\includegraphics[width=0.45\textwidth]{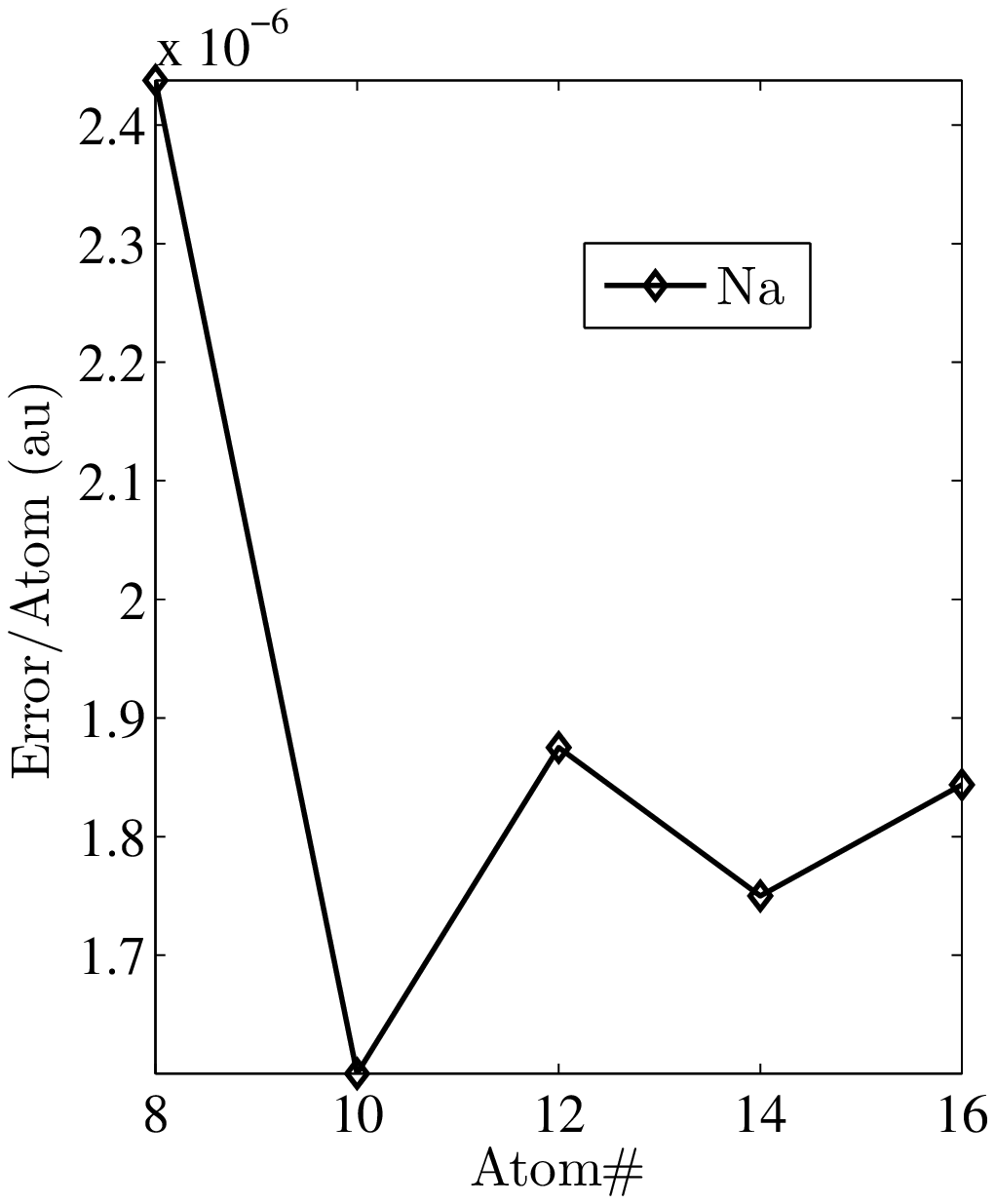}}
    \subfloat[Si]{\includegraphics[width=0.45\textwidth]{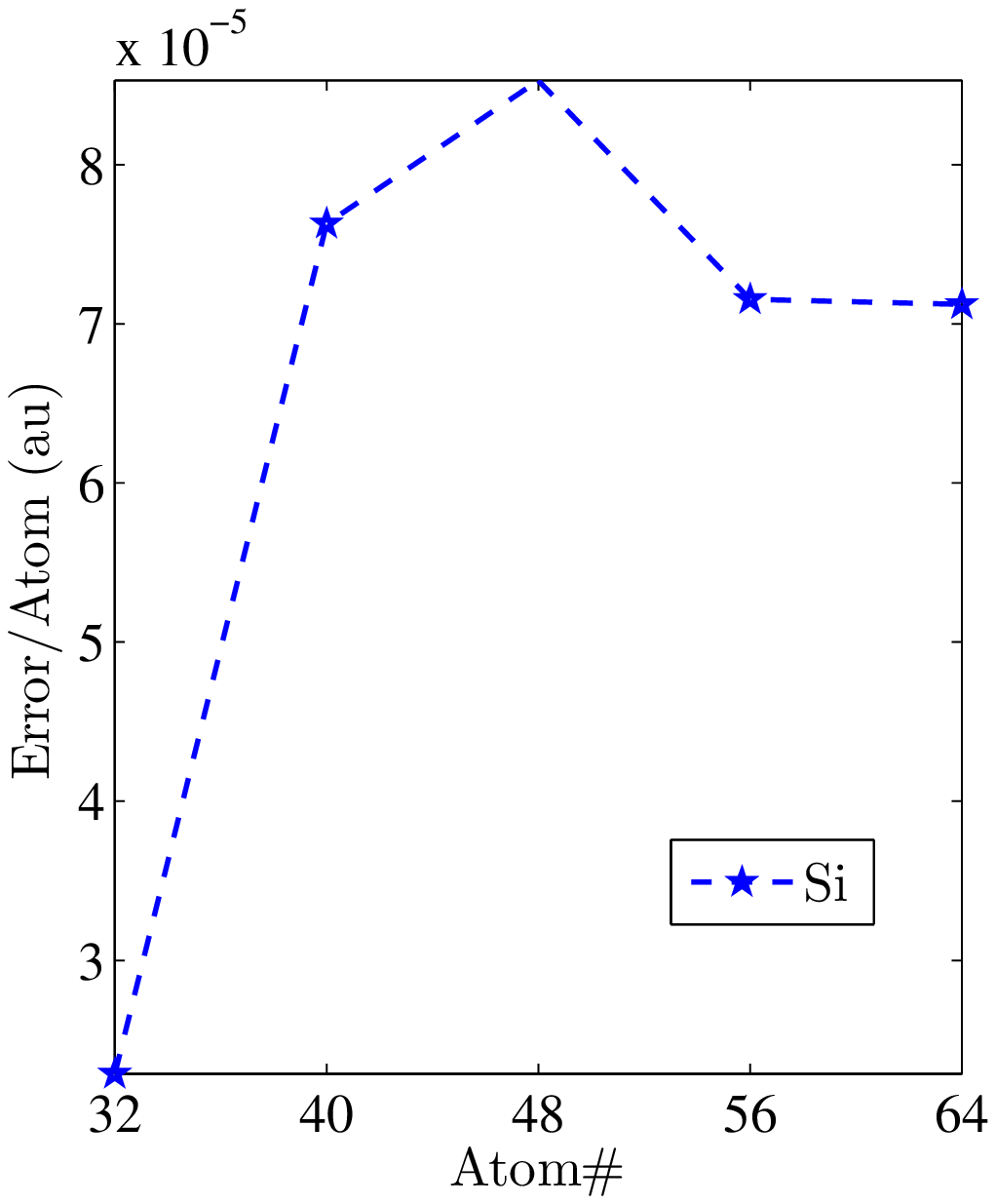}}
  \end{center}
  \caption{(color online) (a) The error of the total energy per atom
    (the $y$ axis) for disordered quasi-1D sodium systems of different
    numbers of atoms (the $x$ axis) extended along the $z$ direction.
    The buffer size is $0.50$, and $4$ adaptive local basis functions
    per atom are used in each calculation.  (b) The error of the total
    energy per atom for the disordered quasi-1D silicon systems of
    different numbers of atoms (the $x$ axis) extended along the $z$
    direction.  The buffer size is $0.50$, and $6$ adaptive local
    basis functions per atom are used in each calculation.}
  \label{fig:NaSi_1Dlong} 
\end{figure}

The behavior of the error found above depends weakly on the number of
atoms of the quasi-1D system extended along the $z$ direction.  The
error of the total energy per atom for disordered quasi-1D systems of
different numbers of atoms is shown for Na in
Fig.~\ref{fig:NaSi_1Dlong} (a) and for Si in
Fig.~\ref{fig:NaSi_1Dlong} (b), respectively.  In both cases the
buffer size is $0.50$, and the penalty parameter is $20$.  Here $4$
and $6$ adaptive local basis functions per atom are used for Na and
Si, respectively.

\subsection{Disordered Quasi-2D and 3D Bulk Systems}\label{subsec:2D3D}

\begin{figure}[ht]
  \begin{center}
    \includegraphics[width=0.6\textwidth]{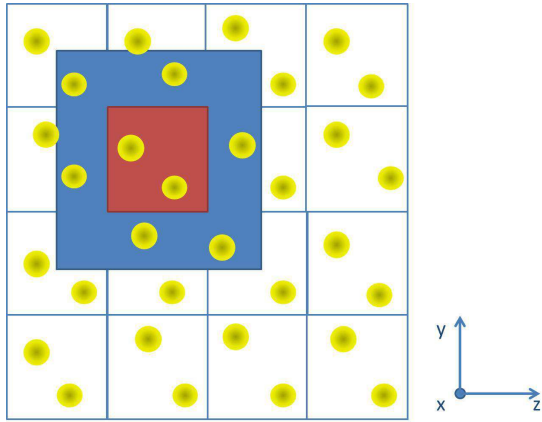}
  \end{center}
  \caption{(color online) A quasi-2D disordered Na system with 32 atoms
  extended along the $y$ and the $z$ directions, viewed along the $x$
  direction.  The red area represents one of the
  elements $E_{2,2}$, and the corresponding extended element $Q_{2,2}$ consists of
  both the red area and the blue area (buffer). The buffer size is $0.50$
  (in the unit of lattice constant) along the $y$ and the $z$ directions, and is
  $0.0$ along the $x$ direction.} 
  \label{fig:Na2Dpartition}
\end{figure}

\begin{figure}[ht]
  \begin{center}
    \subfloat[Quasi-2D Na]{\includegraphics[width=0.45\textwidth]{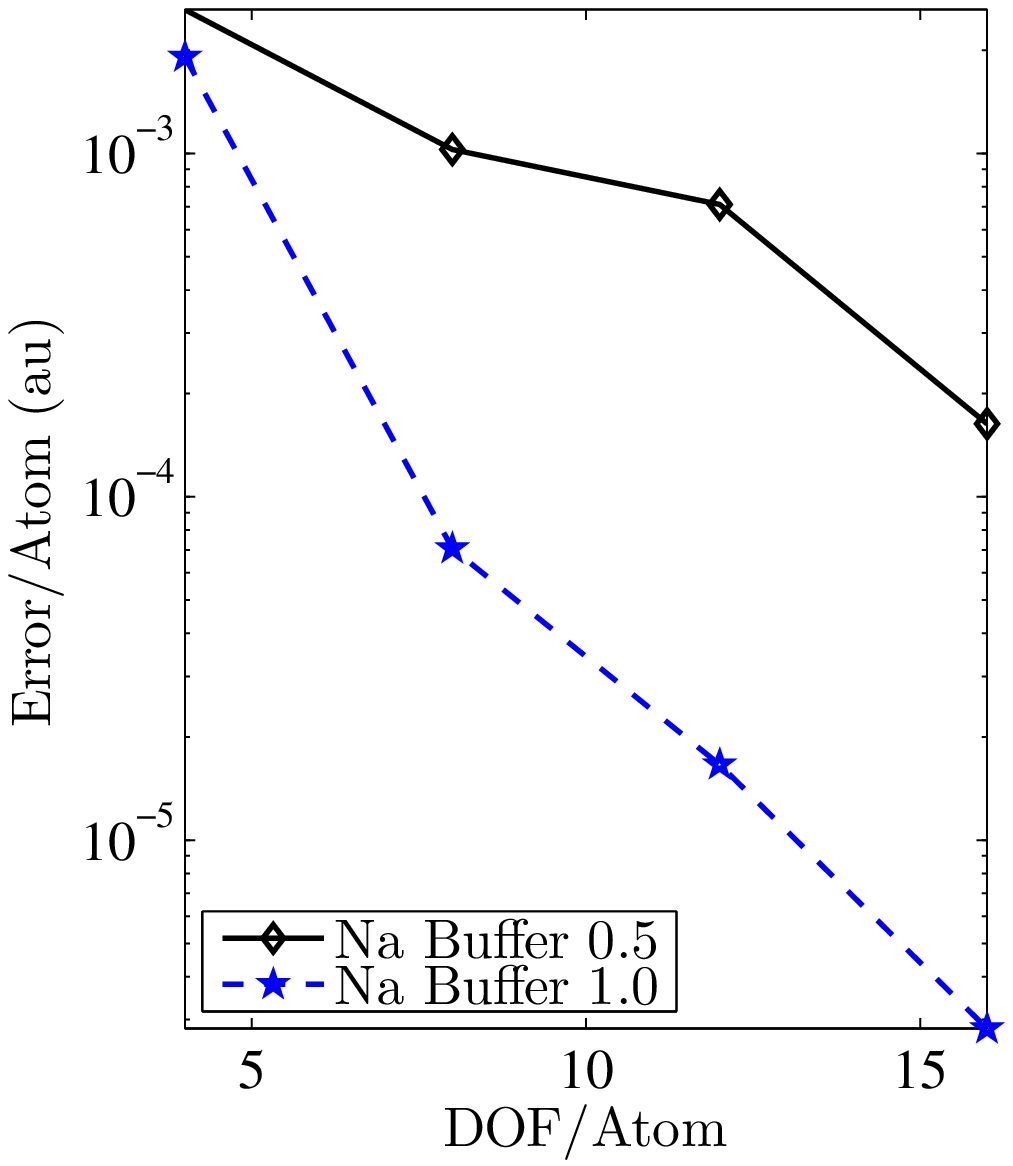}}
    \subfloat[bulk 3D Na]{\includegraphics[width=0.45\textwidth]{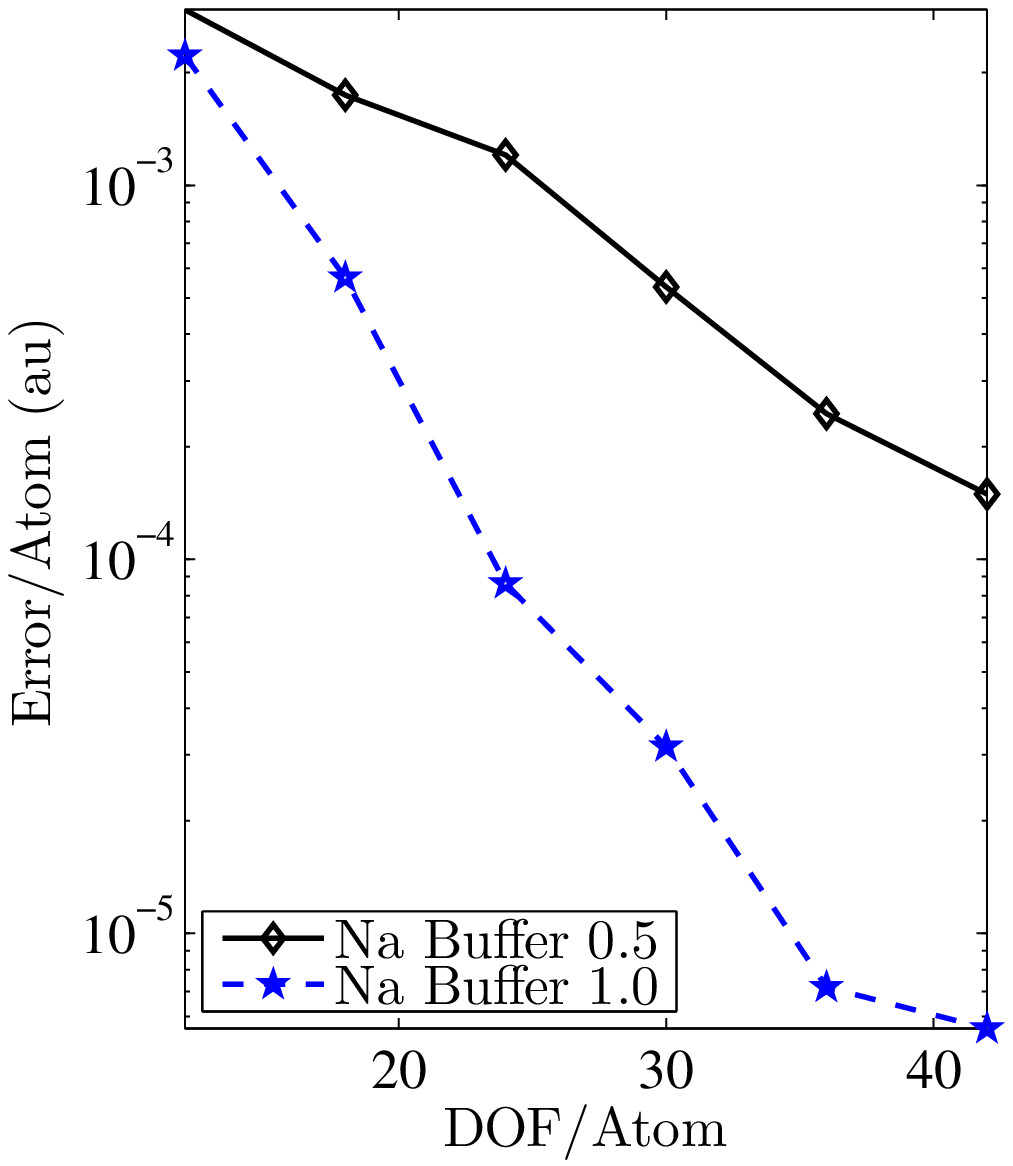}}
  \end{center}
  \caption{(color online) (a) The error of the total energy per atom
  (the $y$ axis, plotted in log-scale) for a disordered quasi-2D sodium system containing $32$
  atoms, with respect to the number of
  basis functions per atom (the $x$ axis).  The  buffer size is chosen
  to be $0.50$ (black diamond with solid line), and $1.00$ (blue star with dashed line),
  respectively.
  (b)  The error of the total energy per atom for a disordered bulk 3D sodium system (the
  $y$ axis, plotted in log-scale) containing $128$ atoms, with respect to the number of
  basis functions per atom (the $x$ axis).
  The buffer size is chosen to be $0.50$ (black diamond with solid line),
  and $1.00$ (blue star with dashed line), respectively. }
  \label{fig:Na_2D3D} 
\end{figure}

This section studies the relation between the error of the total energy
per atom and the dimensionality of the system. The
partition of the domain for systems of higher dimension is similar to
that in the quasi-1D case.
Fig.~\ref{fig:Na2Dpartition} shows the partition of a quasi-2D system with
$32$ sodium atoms, viewed along the $x$ direction.  The domain is partitioned
into $16$ disjoint elements.  The length of each element
(red area) is equal to the length of the lattice constant of the
crystalline unit cell.  The corresponding extended element for solving
the adaptive local basis functions includes both the element (red area) and the
buffer (blue area).  
Fig.~\ref{fig:Na_2D3D} (a) shows the behavior of the error for a
disordered quasi-2D sodium system containing $32$ atoms with the buffer
of size $0.50$ (black diamond with solid line) and of size $1.00$ (blue star with dashed line),
respectively.  For the case with the buffer size equal to $0.50$, the error of the total
energy per atom is $1.0\times 10^{-3}$ au using $8$ basis functions per
atom.  The error of the total energy per atom can reach $2.8\times
10^{-6}$ au with $16$ basis functions per atom and buffer size $1.00$.
Fig.~\ref{fig:Na_2D3D} (b) shows the behavior of the
error for a disordered bulk 3D sodium system containing $128$ atoms with
the buffer of size $0.50$ (black diamond with solid line) and of size $1.00$ (blue star with dashed line),
respectively.
For the case with the buffer size equal to  $0.50$, the error of the total
energy per atom is $1.2\times 10^{-3}$ au using $24$ basis functions per
atom.  The error of the total energy per atom can reach $5.6\times
10^{-6}$ au or $0.15$ meV with $42$ basis functions per atom and buffer size $1.00$.
Compared to the quasi-1D case, the number of adaptive local basis
functions per atom increases significantly in order to reach the same accuracy.  
The increasing number of basis functions is partly due to the increasing
number of Na atoms in the extended
element.  In this case,
the numbers of the Na atoms in the extended element with a buffer
size of $1.00$ are $4,18,54$ for quasi-1D, quasi-2D and bulk 3D
systems, respectively.  The increased number of Na atoms in the extended
elements requires more eigenfunctions in the extended elements, and
therefore more adaptive local basis functions per atom in the elements.  

\subsection{The penalty parameter}\label{subsec:penalty}

\begin{figure}[ht]
  \begin{center}
    \subfloat[Na]{\includegraphics[width=0.45\textwidth]{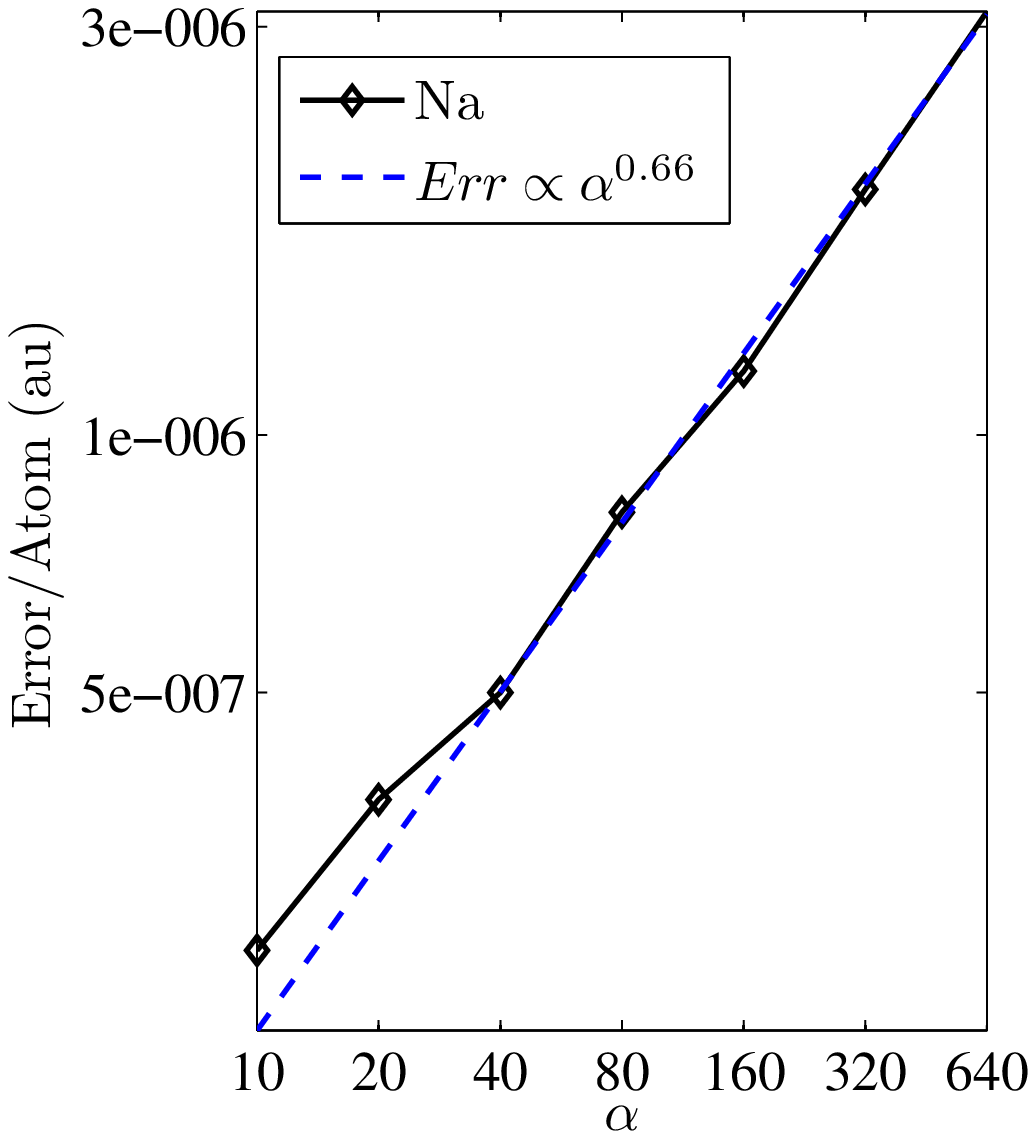}}
    \subfloat[Si]{\includegraphics[width=0.45\textwidth]{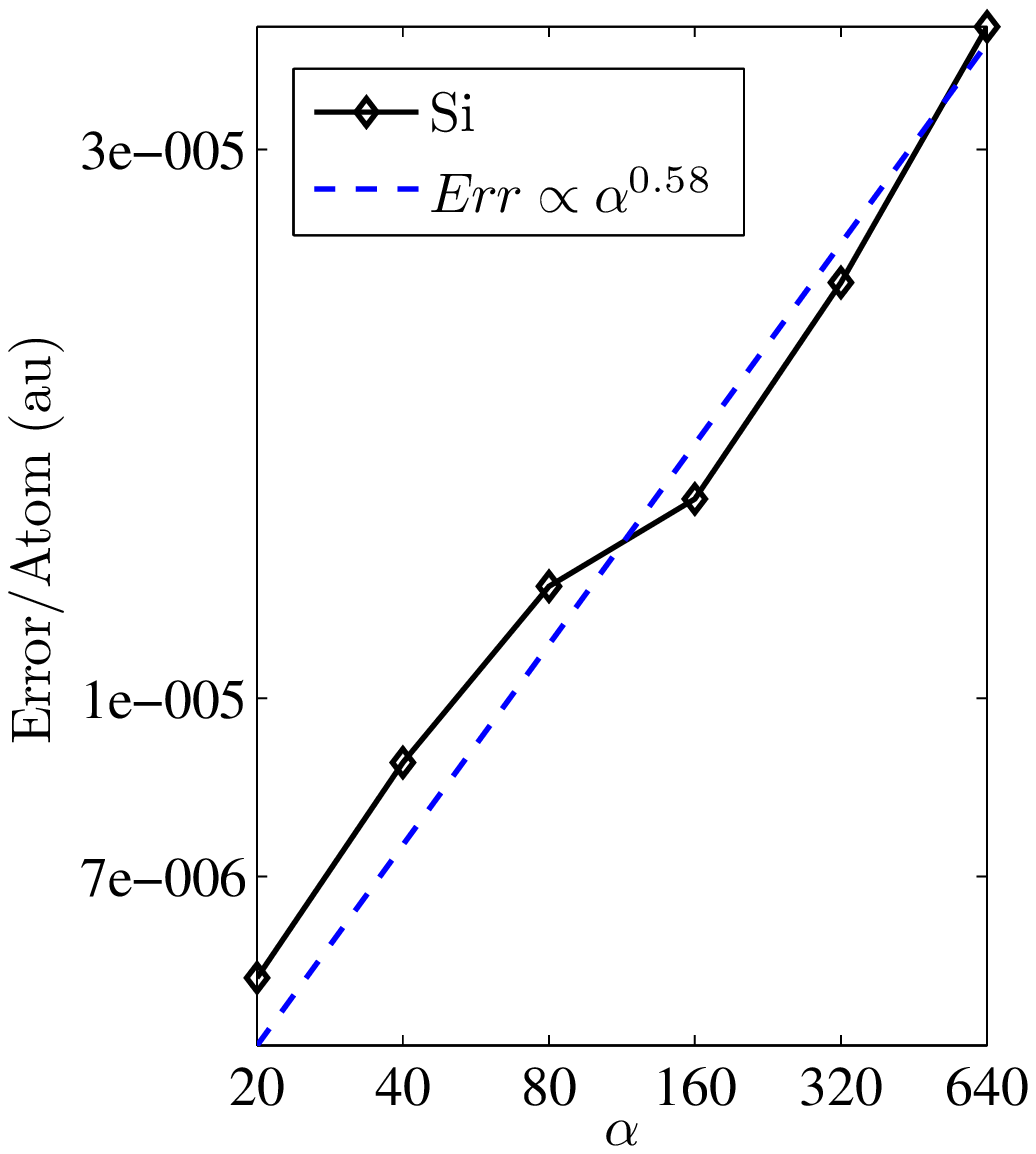}}
  \end{center}
  \caption{(color online) (a) Log-log plot for the error of the total energy per atom
  (the $y$ axis) with respect to the penalty parameter $\alpha$ (the $x$
  axis), for a quasi-1D sodium system with $8$ atoms. The buffer size is
  $1.00$ and the number of basis functions per atom is $12$. The error
  (black diamond with solid line) can be fitted with a polynomial
  function of $\alpha$ (blue dashed line).  (b) Log-log plot for the error of the total
  energy per atom (the $y$ axis) with respect to the penalty parameter
  $\alpha$ (the $x$ axis), for a quasi-1D silicon system with $32$
  atoms. The buffer size is $1.00$ and the number of basis functions per
  atom is $6$. The error (black diamond with solid line) can be fitted
  with a polynomial function of $\alpha$ (blue dashed line).} 
  \label{fig:NaSi_alpha} 
\end{figure}

The interior penalty formulation of the discontinuous Galerkin method
contains an important parameter $\alpha$ to guarantee stability.
$\alpha=20$ has been applied uniformly to all the examples studied
so far.  The $\alpha$-dependence of the
error of the total energy per atom is shown for the quasi-1D sodium
system in Fig.~\ref{fig:NaSi_alpha} (a) and for the quasi-1D silicon
system in Fig.~\ref{fig:NaSi_alpha} (b), respectively.  
For Na, the buffer size is $1.00$, and the number of basis
functions per atom is $8$.  The error of the total energy per atom
is empirically proportional to $\alpha^{0.66}$ up to $\alpha=640$.  For
Si, the buffer size is $1.00$, and the number of basis functions per
atom is $6$. The error of the total energy per atom
is empirically proportional to $\alpha^{0.58}$ up to $\alpha=640$.  
We also remark that the DG formulation can become unstable for 
$\alpha$ smaller than a certain threshold value.  For example, the error
of the total energy per atom is $2.9\times 10^{-1}$ au for Na with
$\alpha=5$, and the error of the total energy per atom is $1.7\times
10^{-2}$ au for Si with $\alpha=10$.  Therefore the penalty parameter
$\alpha$ plays an important role in the stability of the algorithm, but
the DG scheme can be accurate and stable with respect to a large range of
$\alpha$ values.

\subsection{Computational efficiency}\label{subsec:efficiency}

\begin{figure}[ht]
  \begin{center}
    \includegraphics[width=0.60\textwidth]{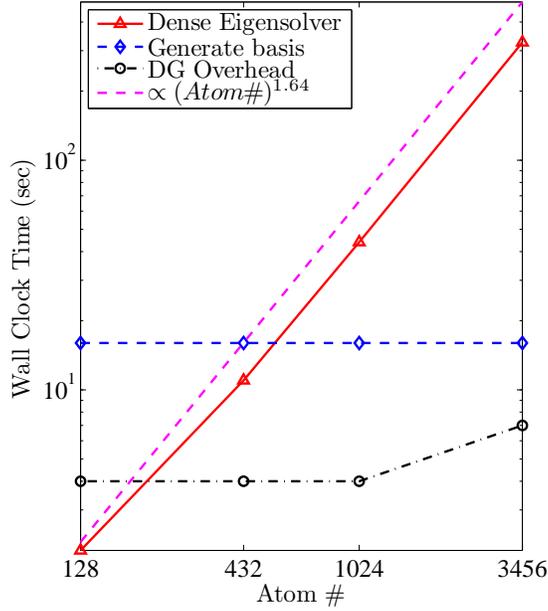}
  \end{center}
  \caption{(color online) Log-log plot for the wall clock time ($y$ axis) for solving
  disordered bulk 3D sodium systems of different sizes ($x$ axis) with
  one step self-consistent field iteration. The number of processors is
  chosen to be proportional to the number of atoms, with $1,728$
  processors used for the largest problem solved here ($3,456$ Na atoms).
  The total wall clock time is broken down into the time for solving the
  DG eigenvalue problem using ScaLAPACK function \textsf{pdsyevd} (red triangle
  with solid line),  the time for generating the adaptive local basis
  functions in the extended elements using LOBPCG solver (blue diamond
  with dashed line), and the time for the overhead in the DG
  calculation, including the matrix assembly and data communication
  (black circle with dot dashed line).  The buffer size is $1.00$, and
  the number of basis functions per atom is $16$.  The scaling of the
  wall clock time for solving the DG eigenvalue problem using \textsf{pdsyevd}
  with respect to the number of atoms is illustrated by the magenta
  dashed line.}
  \label{fig:Na3DTime} 
\end{figure}

The small number of the adaptive basis functions per atom can lead to
significant savings of the computational time.  We illustrate the
efficiency of our algorithm using a disordered bulk 3D sodium
system with the buffer size of $1.00$ and with $16$ basis functions
per atom.  Fig.~\ref{fig:Na_2D3D} (b) suggests that 
the error of the total energy per atom is about $10^{-3}$ au for
this choice of the parameters.  The size of each
element is equal to the lattice constant with $2$ Na atoms in
each element. The size of the global domain
$\Omega$ ranges from $4\times4\times4$ unit cells with $128$ Na atoms
to $12\times12\times12$ elements with $3,456$ atoms.  The number of
processors (cores) used is
proportional to the number of elements, and $1,728$ processors are
used in the problem with $12\times12\times12$ elements.  We compare
the wall clock time for one step self consistent iteration with $3$
LOBPCG iterations for solving the adaptive basis functions in the
extended elements.  Fig.~\ref{fig:Na3DTime} compares the wall clock time
for solving the DG eigenvalue problem using ScaLAPACK function
\textsf{pdsyevd} (red triangle
with solid line),  the time for generating the adaptive local basis
functions in the extended elements using LOBPCG solver (blue diamond
with dashed line), and the time for the overhead in the DG
calculation (black circle with dot dashed line).  The buffer size is $1.00$, and
the number of basis functions per atom is $16$.  
Since both the size of the extended elements and the number of basis
functions per atom are fixed, the computational time for solving the
adaptive basis functions does not depend on the global domain size.
The overhead in the DG calculation method includes mainly the assembly
process of the DG Hamiltonian matrix via numerical integration and data
communication.  All numerical integrations are
localized inside each element and its neighboring elements.  Our
implementation ensures that the data communication is restricted to be
within nearest neighboring elements.  Therefore the time for the
overhead increases mildly with respect to the global system size. 
The complexity of the DG eigensolver using \textsf{pdsyevd} scales
cubically with respect to global system size in the asymptotic limit,
and starts to dominate the cost of computational time for system
containing more than $1,000$ atoms.  
Since the number of processors is proportional to the number
of elements, the asymptotic wall clock time for the DG eigensolver
should scales quadratically with respect to the number of atoms.  
The practical wall clock time for solving the DG eigensolver is found to
be proportional  to $(N_{atom})^{1.64}$ (magenta dashed line in
Fig.~\ref{fig:Na3DTime}), indicating that the asymptotic cubic scaling
has not yet been reached.  In the largest example with $3,456$ atoms,  the matrix size of
the DG Hamiltonian matrix is $55,296$.

The efficiency due to the dimension reduction of the adaptive basis
functions can be illustrated by the comparison between the cost of the
computational time of the LOBPCG eigensolver in the global
domain calculation (Global), and that of the DG
eigenvalue problem with the adaptive basis functions (DG), as reported
in Table~\ref{tab:ctime}.  The global domain calculation uses 
$10$ LOBPCG iteration steps per SCF iteration.  On a single processor,
the global domain calculation 
costs $806$ sec for the bulk 3D sodium system with $128$ atoms, and
$19,112$ sec for the bulk 3D sodium system with $432$ atoms.  By
assuming that the global domain calculation can be ideally parallelized, the third
column of Table~\ref{tab:ctime} reports the computational time of the
global domain calculation
measured on a single processor divided by the number of processors
used in the corresponding DG eigensolver.  The fourth column reports
the wall clock time for the DG eigensolver executed in parallel. We
remark that the computational time for solving the adaptive local
basis functions is not taken into account, since we are comparing the
savings of the computational time due to the dimension reduction of
the basis functions.  It is found that the saving of the computational
time is already significant even when the system size is relatively
small.

\begin{table}[h]
  \centering
  \begin{tabular}{c|c|c|c}
    \hline
    Atom\# & Proc\# & Global (sec) & DG (sec)\\
    \hline
    128 & 64  & 13 &  1\\
    432 & 216 & 88 & 14\\
    \hline
  \end{tabular}
  \caption{The comparison of the cost of the wall clock time using
  the LOBPCG iteration on the global domain (performed with a single
  processor and divide the time by the number of processors in column
  $2$, assuming that the LOBPCG are perfectly parallelized) and the wall
  clock time using the adaptive local basis functions (only count the DG
  eigenvalue solver using ScaLAPACK with the number of processors in
  column $2$). The systems under study are the bulk 3D sodium system
  with $4\times4\times4$ elements ($128$ Na atoms), and with
  $6\times6\times6$ elements ($432$ Na atoms), respectively. 
  }
  \label{tab:ctime}
\end{table}

\section{Discussion and Conclusion}\label{sec:discussion}

In this paper we proposed the adaptive local basis functions for
discretizing the Kohn-Sham Hamiltonian operator, and demonstrated that
the adaptive local basis functions are efficient for calculating the
total energy and electron density, and can reach high accuracy
with a very small number of basis functions per atom.  The adaptive
local basis functions are discontinuous in the global domain, and the
continuous Kohn-Sham orbitals and electron density are reconstructed
from these discontinuous basis functions using the discontinuous
Galerkin (DG) framework.  The environmental effect is automatically
built into the basis functions, thanks to the flexibility provided by
the DG framework.

The current implementation of the DG method is already able to perform
the total energy calculation for systems consisting of thousands of
atoms.  The performance of the DG method can be improved by taking
into account the block sparsity of the DG stiffness matrix.
Furthermore, the local nature of the adaptive basis functions allows
us to incorporate the recently developed pole expansion and selected
inversion type fast algorithms~\cite{LinLuYingE2009,LinLuYingEtAl2009,LinYangLuEtAl2011,LinYangMezaEtAl2010} into the DG framework.  The
capability of the resulting algorithm is expected to be greatly
enhanced compared to the current implementation.  This is our ongoing
work.

In order to generalize the current framework to the force calculation
and further to the geometry optimization and the \textit{ab initio}
molecular dynamics simulation, the adaptive local basis functions and
their derivatives with respect to the positions of the atoms (called
Pulay force~\cite{Pulay:69}) should be both accessible.  Recently we
propose the optimized local basis functions~\cite{LinLuYingE2011_opt}
that is able to systematically control the magnitude of the Pulay
force, which is a further improvement of the adaptive local basis
functions.  This is also our ongoing work.


\vspace{1em}
\noindent{\bf Acknowledgement:}

This work is partially supported by DOE under Contract No.
DE-FG02-03ER25587 and by ONR under Contract No. N00014-01-1-0674 (W.~E
and L.~L.), and by a Sloan Research Fellowship and by NSF CAREER grant
DMS-0846501 (L.~Y.).  We thank the National Energy Research Scientific
Computing Center (NERSC), and the Texas Advanced Computing Center (TACC)
for the support to perform the calculations. L.~L. and J.~L. thank the
University of Texas at Austin for the hospitality where the idea of this
paper starts.



\FloatBarrier


\end{document}